\documentclass[12pt]{article}
\usepackage[final]{epsfig}
\usepackage{graphics}
\usepackage{amsmath}
\usepackage{amsfonts}
\usepackage{latexsym}
\usepackage{amssymb}
\usepackage{graphicx}
\usepackage{url}
\usepackage{epstopdf}

\newtheorem{lemma}{Lemma}[section]
\newtheorem{proposition}[lemma]{Proposition}
\newtheorem{remark}[lemma]{Remark}

\newtheorem{theorem}{Theorem}[section]

\begin{document}
\newcommand{\eps}{{\varepsilon}}
\newcommand{\proofend}{$\Box$\bigskip}
\newcommand{\C}{{\mathbb C}}
\newcommand{\Q}{{\mathbb Q}}
\newcommand{\R}{{\mathbb R}}
\newcommand{\Z}{{\mathbb Z}}
\newcommand{\RP}{{\mathbb {RP}}}
\newcommand{\CP}{{\mathbb {CP}}}
\newcommand{\Tr}{\rm Tr}
\def\proof{\paragraph{Proof.}}

\title{Projective configuration theorems: old wine into new wineskins}

\author{Serge Tabachnikov\footnote{
Department of Mathematics,
Penn State University,
University Park, PA 16802;
tabachni@math.psu.edu}
}

\date{}
\maketitle

\tableofcontents

\section{Introduction: classical configuration theorems} \label{intro}

Projective configuration theorems are among the oldest and best known mathematical results. The next figures depict the famous theorems of Pappus, Desargues, Pascal, Brianchon, and Poncelet.

\begin{figure}[hbtp]
\centering
\includegraphics[height=1.7in]{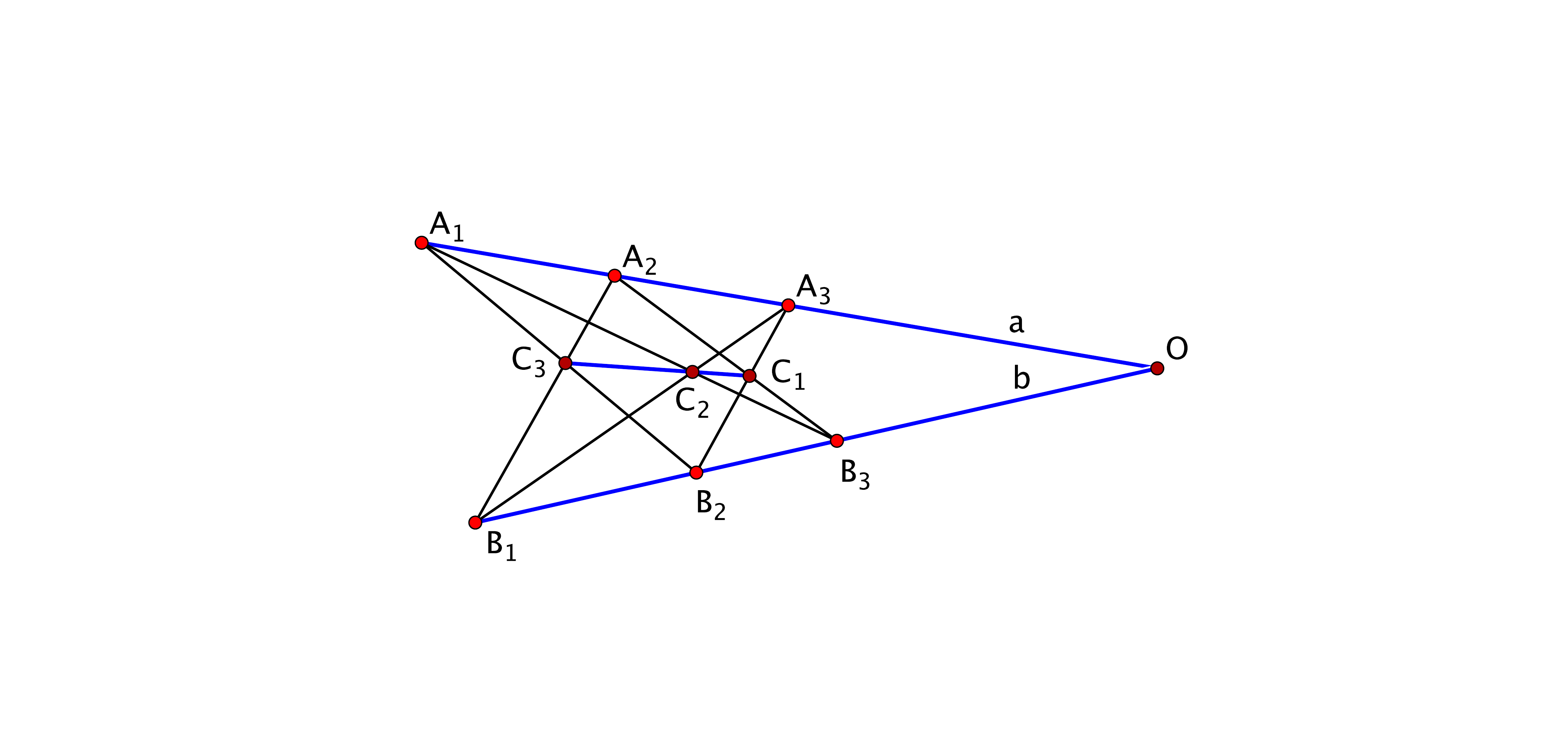}
\caption{The Pappus theorem: if $A_1,A_2,A_3$ and $B_1,B_2,B_3$ are two collinear triples of points, then $C_1,C_2,C_3$ is also a collinear triple.}
\label{Pappusfig}
\end{figure}

\begin{figure}[hbtp]
\centering
\includegraphics[height=2.6in]{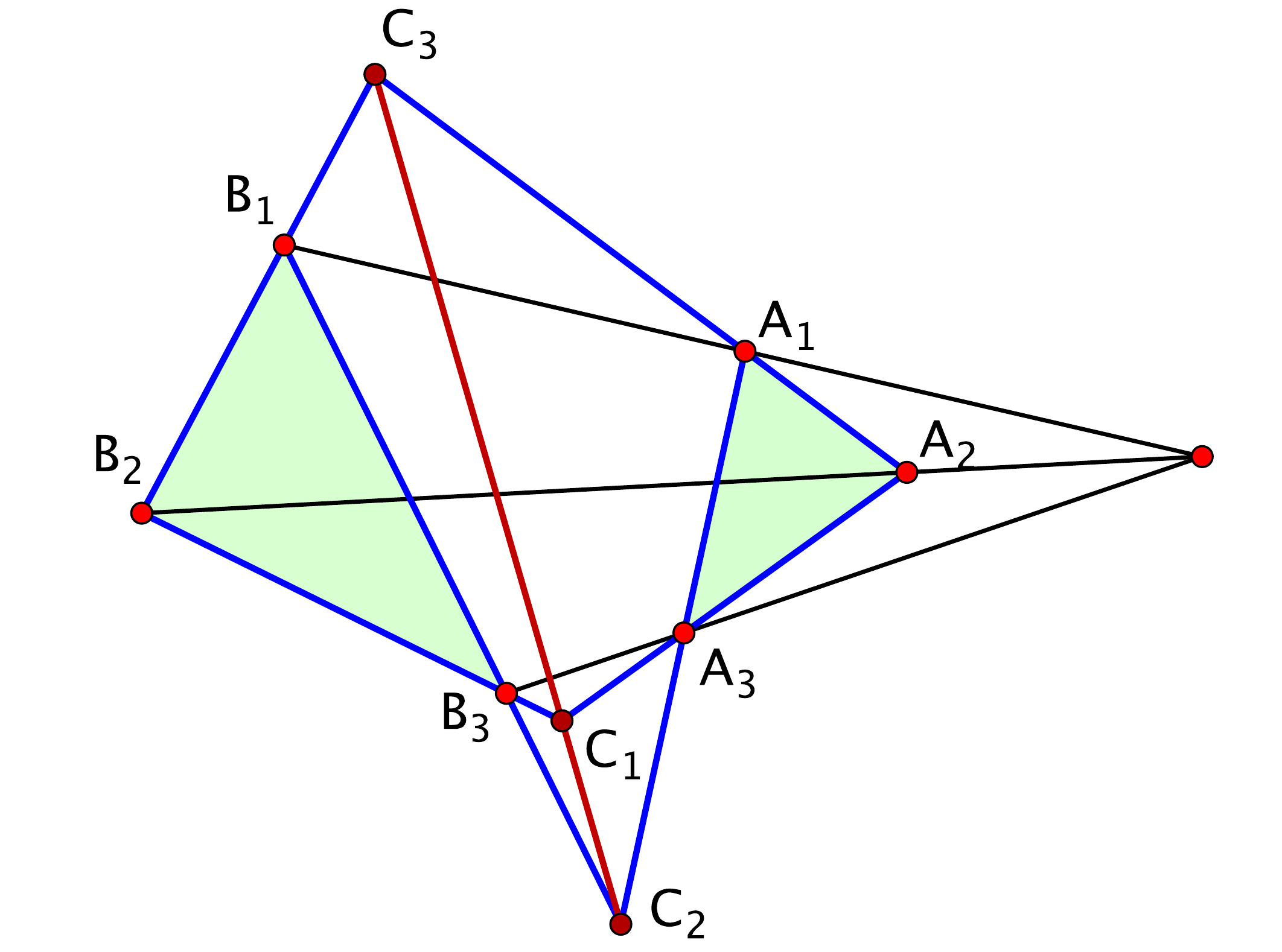}
\caption{The Desargues theorem: if the lines $A_1 B_1, A_2 B_2$ and $A_3 B_3$ are concurrent, then the points  $C_1,C_2,C_3$ are collinear.}
\label{Desarguesfig}
\end{figure}

\begin{figure}[hbtp]
\centering
\includegraphics[width=2.8in]{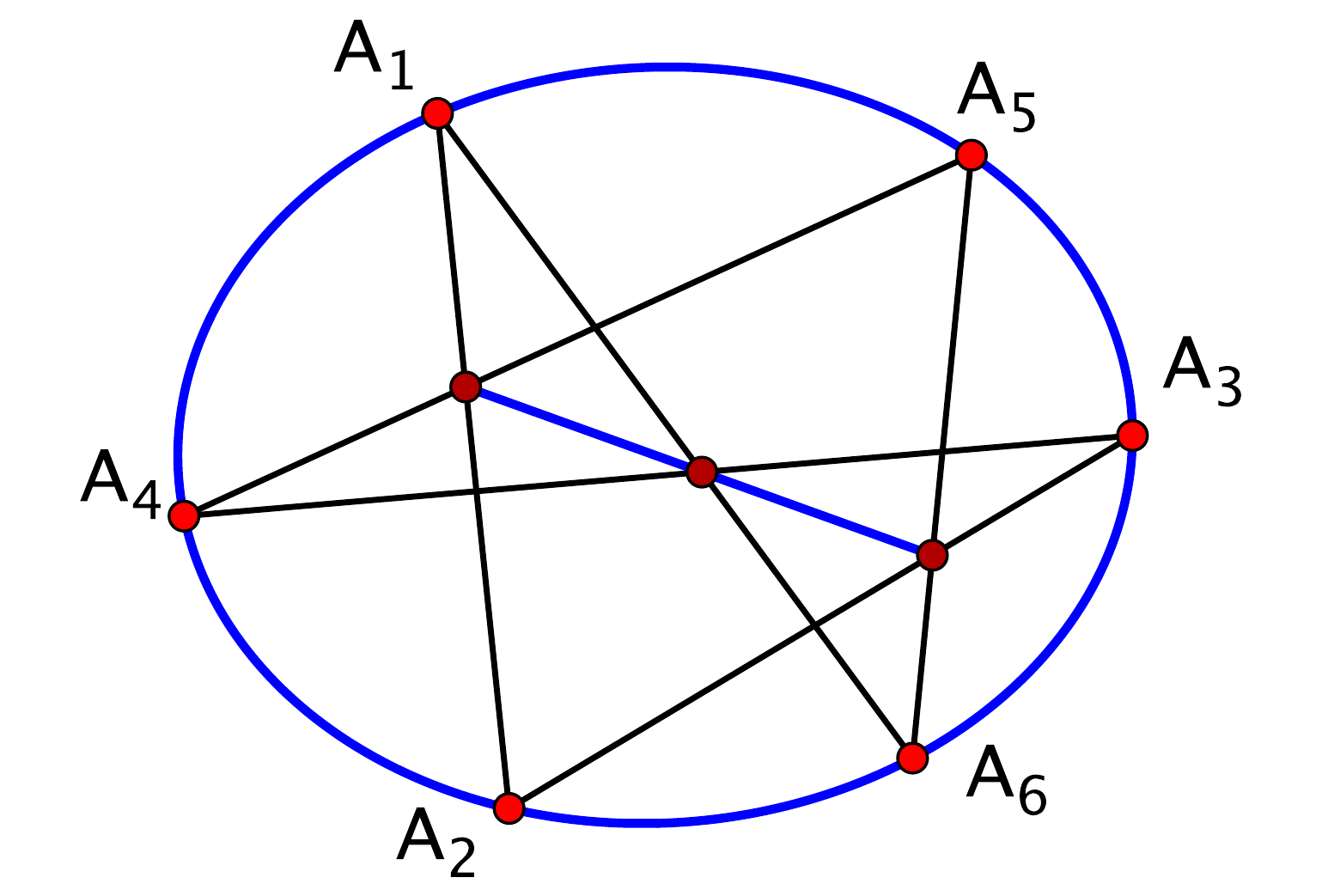} \
\includegraphics[width=2.5in]{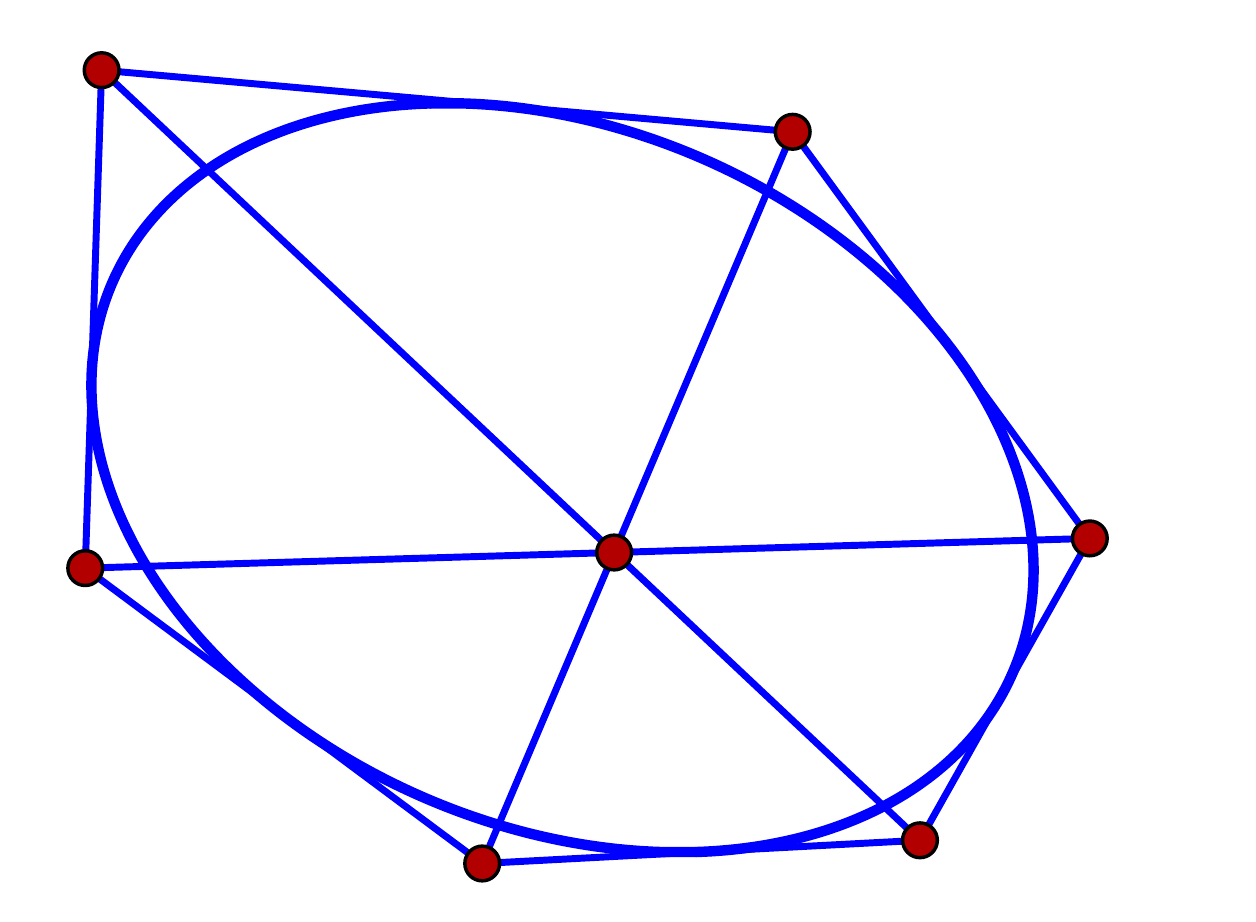}
\caption{The Pascal theorem, a generalization of the Pappus theorem: the points $A_1, \ldots, A_6$  lie on a conic, rather than the union of two lines. The Brianchon theorem is projectively dual to Pascal's.}
\label{Pascalfig}
\end{figure}

\begin{figure}[hbtp]
\centering
\includegraphics[height=3in]{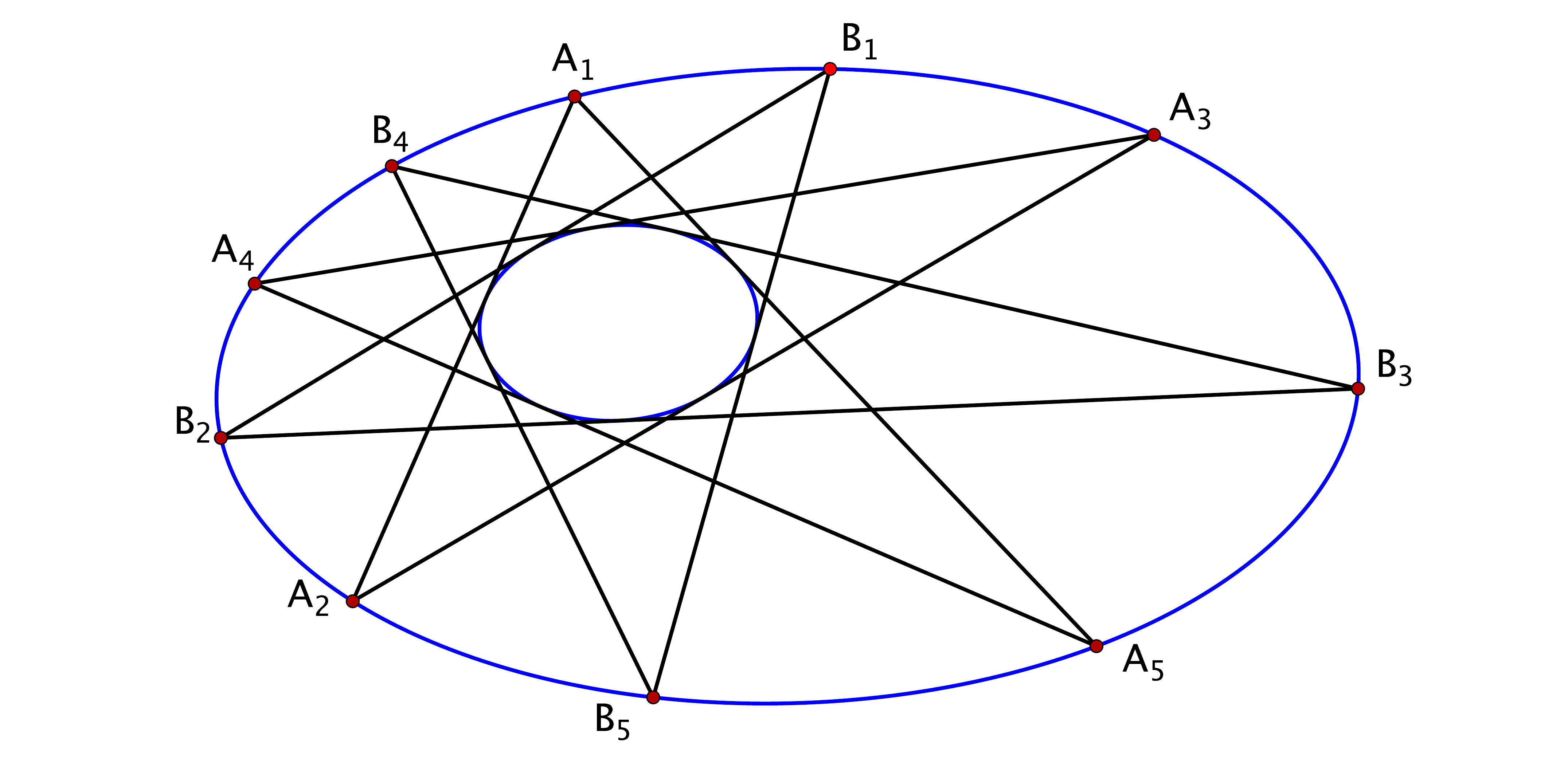}
\caption{The Poncelet Porism, case $n=5$: if the polygonal line $A_1 A_2 A_3 A_4 A_5$, inscribed into a conic and circumscribed about a conic, closes up after five steps, then so does any other polygonal line $B_1 B_2 B_3 B_4 B_5$.}
\label{Ponceletfig}
\end{figure}

The literature on configuration theorems is vast; the reader interested in a panoramic view of the subject is recommended \cite{Be,RG}. 

Configuration theorems continue to be an area of active research.  
To a great extent, this is due to the advent of computer as a tool of experimental research in mathematics. In particular, interactive geometry software is a convenient tool for the study of geometric configurations. The illustration in this article are made using such a software, Cinderella 2 \cite{Cin}.

Another reason for the popularity of configuration theorems is that they play an important role in the emerging field of discrete differential geometry and the theory of completely integrable systems \cite{BS}.

The goal of this survey is to present some recent results motivated and inspired by the classical configuration theorems;  these results make the old theorems fresh again.
The selection of topics reflects this author's taste; no attempt was made to present a  comprehensive description  of the area. In the  cases when proofs are discussed, they are only outlined; the reader interested in details is referred to the original papers.

We assume that the reader is familiar with the basics of projective, Euclidean, spherical, and hyperbolic geometries. One of the standard references is \cite{Be1}, and \cite{HC} is as indispensable as ever.

Now let us specify what we mean by configuration theorems in this article. 
The point of view is dynamic, well adapted for using interactive geometry software.

An initial data for a configuration theorem is a collection of labelled points $A_i$ and lines $b_j$ in the projective plane, such that, for some pairs of indices $(i,j)$, one has the incidence $A_i \in b_j$. If, in addition, a polarity is given, then one can associate the dual line to a point, and the dual point to a line. In presence of polarity, the initial data includes information that, for some pairs of indices $(k,l)$,  the point $A_k$ is polar dual to the line $b_l$.

One also has an ordered list of instructions consisting of two operations: draw a line through a  pair of  points, or intersect a  pair of  lines at a point. These new  lines and points also receive labels. If polarity is involved, one also has the operation of taking the polar dual object, point to line, or line to point.

The statement of a configuration theorem is that, among so constructed points and lines,  certain incidence relations hold, that is, certain points lie on certain lines. 

It is  assumed that the conclusion of a configuration theorem holds for almost every initial set of points and lines satisfying the initial conditions, that is, holds for a Zariski open set of such initial configurations. This is different from what is meant by a configuration of points and lines  in chapter 3 of \cite{HC} or in \cite{Gr}: the focus there is on whether a combinatorial incidence is realizable by  points and lines in the projective plane.

\bigskip
{\bf Acknowledgements}. 
I am grateful to R. Schwartz for numerous stimulating discussions and to P. Hooper for an explanation of his work.
I was supported by NSF grant DMS-1510055. This article was written during my stay at ICERM; it is a pleasure to thank the Institute for its inspiring  and friendly atmosphere. 

\section{Iterated Pappus theorem and the modular group} \label{itPapp}

The Pappus theorem can be viewed as a construction in $\RP^2$ that inputs two ordered triples of collinear points $A_1,A_2,A_3$ and $B_1,B_2,B_3$, and outputs a new collinear triple of points $C_1,C_2,C_3$, see Figure \ref{Pappusfig}. One is tempted to iterate: say, take $A_1,A_2,A_3$ and $C_1,C_2,C_3$ as an input. Alas, this takes one back to the triple $B_1,B_2,B_3$.

\begin{figure}[hbtp]
\centering
\includegraphics[height=2in]{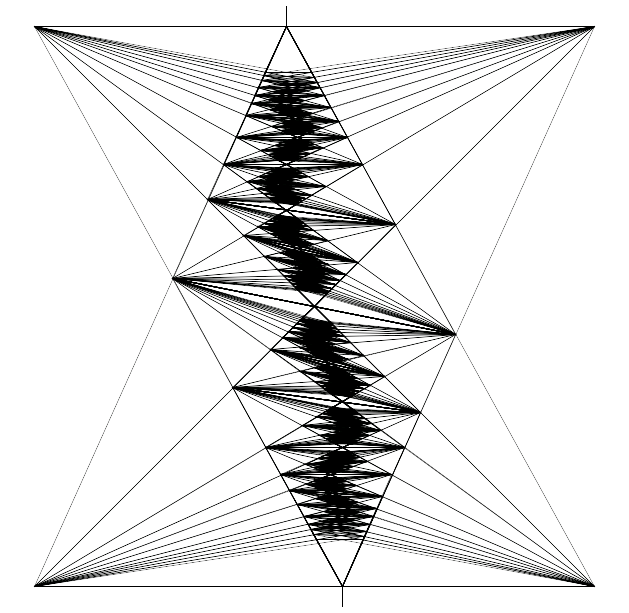}
\caption{Iterations of the Pappus construction produced by R. Schwartz's applet \cite{App}.}
\label{Pappusiter}
\end{figure}

To remedy the situation, swap points $C_1$ and $C_3$. Then the input $A_1,A_2,A_3$ and $C_1,C_2,C_3$ yields a new collinear triple of points, and so does the input $C_1,C_2,C_3$ and $B_1,B_2,B_3$. And one can continue in the same way indefinitely, see  Figure \ref{Pappusiter}. 
The study of these iterations was the topic of R. Schwartz's paper \cite{Sch93}.

Return to Figure \ref{Pappusfig}. The input of the Pappus construction is the {\it marked box} $(A_1,A_3, B_3, B_1; A_2, B_2)$, a quadrilateral $A_1 A_3 B_3 B_1$ with the top distinguished point $A_2$ and the bottom distinguished point $B_2$. The marked box is assumed to satisfy the {\it convexity condition}: the points $A_1$ and $A_3$ are separated by the points $A_2$ and $O$ on the projective line $a$, and likewise for the pairs of points $B_1, B_3$ and $B_2, O$ on the line $b$. Marked boxes that differ by the involution
$$
(A_1, A_3, B_3, B_1; A_2, B_2) \leftrightarrow (A_3, A_1, B_1, B_3; A_2, B_2)
$$
are considered to be the same. 

A convex set in $\RP^2$ is a set that is disjoint from some line and that is convex in the complement to this line, the affine plane.  Two points in $\RP^2$ can be connected by a segment in two ways. 
The four points $A_1, A_3, B_3, B_1$, in this cyclic order, define 16 closed polygonal lines, but only one of them is the boundary of a convex quadrilateral, called the interior of the convex marked box. 

Recall that the points of the dual projective plane  are the lines of the initial plane.
Let $\Theta=(A_1,A_3, B_3, B_1; A_2, B_2)$ be a convex marked box in $\RP^2$. Its  dual, $\Theta^*$, is a marked box in the dual projective plane whose points are the lines
$$
(A_2B_1,A_2B_3,A_1B_2,A_3B_2; a,b).
$$
The dual marked box is also convex.

The moduli space of projective equivalence classes of marked boxes in 2-dimensional. One can send the points $A_1, A_3, B_3, B_1$ to the vertices of a unit square; then the projective equivalence class of a convex marked box is determined by the positions of the points $A_2$ and $B_2$ on the horizontal sides of the square. Namely, let $x=|A_1 A_2|, y=|B_1 B_2|$. Then the equivalence class 
\begin{equation} \label{invol}
(x,y) \sim (1-x,1-y),
\end{equation}
where $0<x,y<1$, determines the the projective equivalence class of a convex marked box. We denote this equivalence class by $[x,y]$.

\begin{figure}[hbtp]
\centering
\includegraphics[height=1.8in]{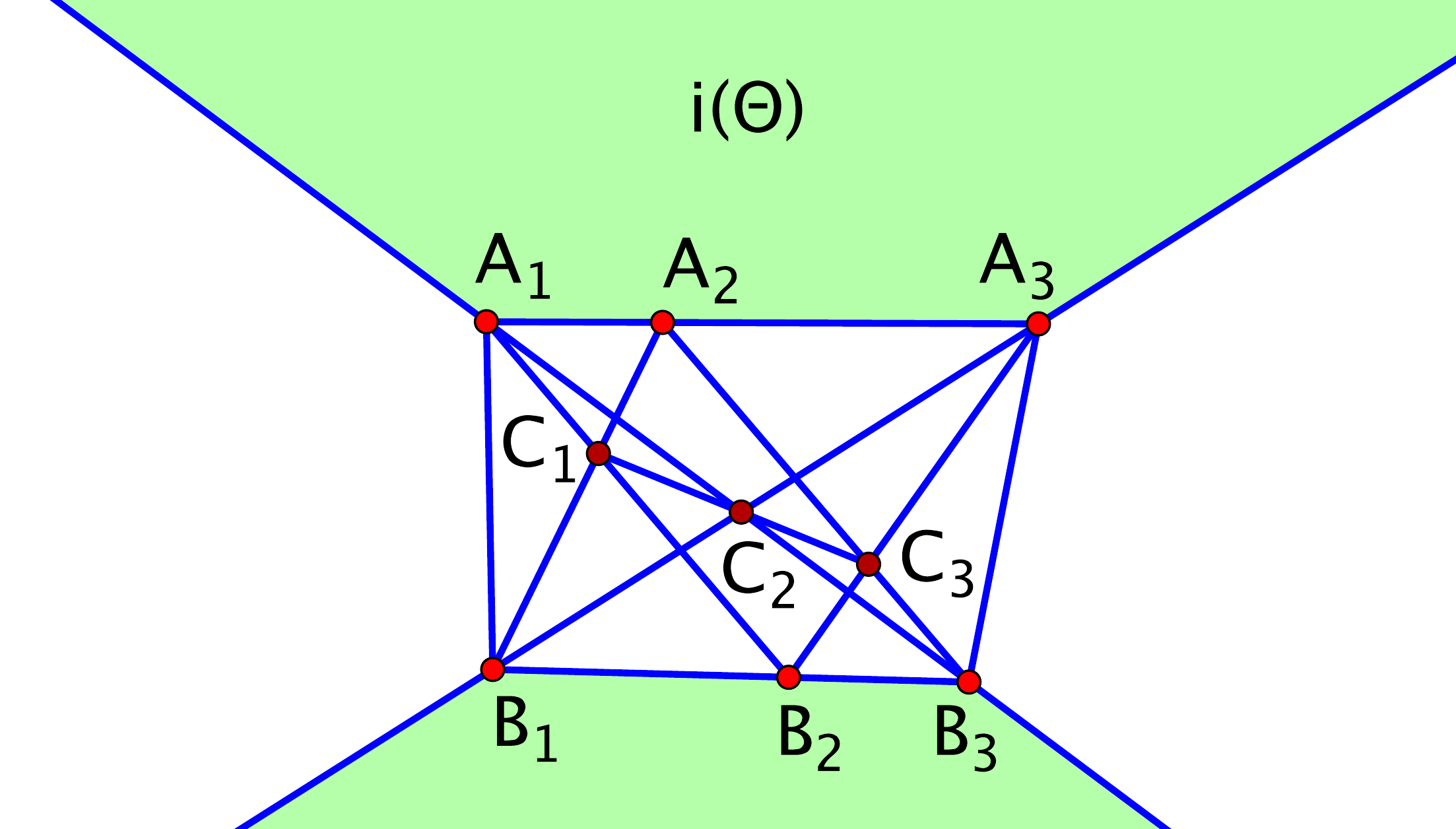}
\caption{The interior of the convex marked box $i(\Theta)$ is bounded by the segments $A_1 A_3, A_3 B_1, B_1 B_3$ and $B_3 A_1$. Two of these segments cross the line at infinity.}
\label{complement}
\end{figure}

The Pappus construction  defines two operations on convex marked boxes, see Figure \ref{complement}:
\begin{equation*}
\begin{split}
\tau_1: (A_1,A_3, B_3, B_1; A_2, B_2) \mapsto (A_1,A_3, C_3, C_1; A_2, C_2),\\
\tau_2: (A_1,A_3, B_3, B_1; A_2, B_2) \mapsto (C_1,C_3, B_3, B_1; C_2, B_2).
\end{split}
\end{equation*}

Add to it a third operation
$$
i: (A_1,A_3, B_3, B_1; A_2, B_2) \mapsto (B_1,B_3,A_1,A_3; B_2,A_2),
$$
also shown in  Figure \ref{complement}.

The three operations form a semigroup $G$. The operations satisfy the following identities,  proved by inspection.

\begin{lemma} \label{groupG}
One has:
$$
i^2=1,\ \tau_1i\tau_2=\tau_2 i \tau_1 = i,\ \tau_1 i \tau_1 = \tau_2,\ \tau_2 i \tau_2 = \tau_1.
$$
\end{lemma}
As a consequence, $G$ is a group; for example, $\tau_1^{-1}= i \tau_2 i$.

Recall that the modular group $M$  is the group of fractional-linear transformations with integral coefficients and determinant one, that is, the group $PSL(2,\Z)$. Since $PGL(2,\R)$ is the group of orientation preserving isometries of the hyperbolic plane, the modular group $M$ 
is realized as a group of isometries of $H^2$. 

\begin{figure}[hbtp]
\centering
\includegraphics[height=2.3in]{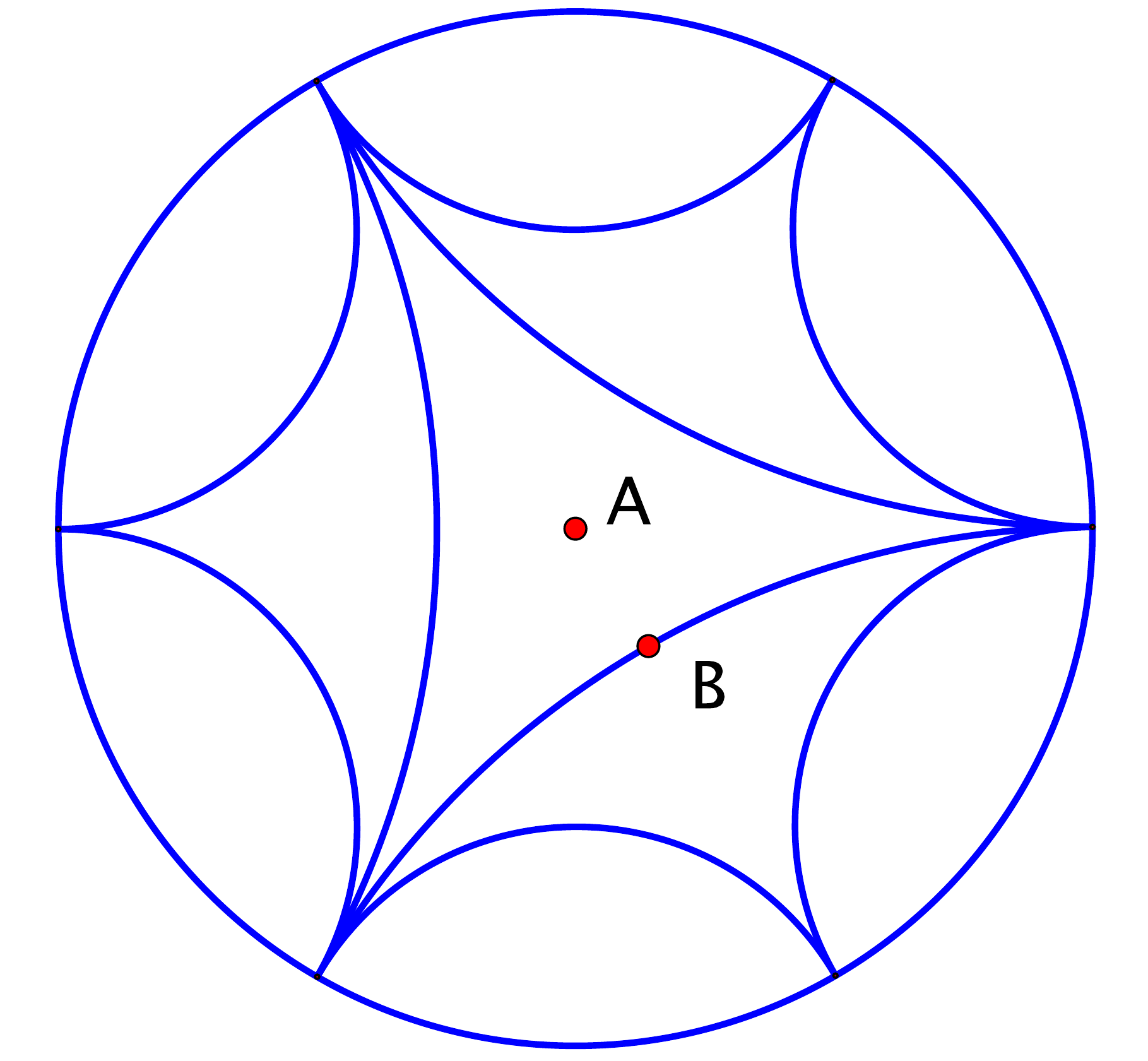}
\caption{A tiling of the hyperbolic plane, in the Poincar\'e disk model, by ideal triangles.}
\label{Farey}
\end{figure}

Consider the tiling of $H^2$ by ideal triangles obtained from one such triangle by consecutive reflections in the sides, see Figure \ref{Farey} for the beginning of this construction. 
The modular group is generated by two symmetries of the tiling: the order three rotation  about point $A$ and the order two rotation (central symmetry)  about point $B$. Algebraically, $M$ is a free product of $\Z_3$ and $\Z_2$.

Return to the group $G$. It is generated by the elements $\alpha=i\tau_1$ and $\beta = i$. Lemma \ref{groupG} implies that $\alpha^3=\beta^2=1$. One can prove that there are no other relations, and hence $G=\Z_3 * \Z_2$ is identified with the modular group.

Given a convex marked box $\Theta$, consider its orbit $\Omega = G(\Theta)$ under the action of the group $G$. The orbit can be described by its oriented incidence graph $\Gamma$. The edges of $\Gamma$ correspond to the marked boxes of $\Omega$, oriented from top to bottom, and the vertices correspond to the tops and the bottoms of the boxes. 

One can embed $\Gamma$ in the hyperbolic plane as in Figure \ref{Farey} (the orientations of the edges are not shown). The group $G$ acts by permutations of the edges of $\Gamma$. The operation $i$ reverses the orientations of the edges. The operation $\tau_1$  rotates each edge counterclockwise one `click' about its tail, and $\tau_2$ rotates the edges one `click' clockwise about their heads. (This is a different action from the one generated by rotations about points $A$ and $B$ in Figure \ref{Farey}). Denote by $G'$ the index two subgroup of $G$ that consists of the transformations that preserve the orientations of the edges.

The orbit $\Omega$ of a convex marked box $\Theta$ has a large group of projective symmetries, namely, an index two subgroup $M'$ of the modular group $M$. This is one of the main results of \cite{Sch93}.
Specifically, one has

\begin{proposition} \label{projsym}
Given a convex marked box $\Theta$, there is an order three projective transformation with the cycle
$$
i(\Theta) \mapsto \tau_1(\Theta) \mapsto \tau_2(\Theta).
$$
In addition, there exists a polarity that identifies $i(\Theta)$ with the  dual  box $\Theta^*$.
\end{proposition}

\begin{figure}[hbtp]
\centering
\includegraphics[height=2.5in]{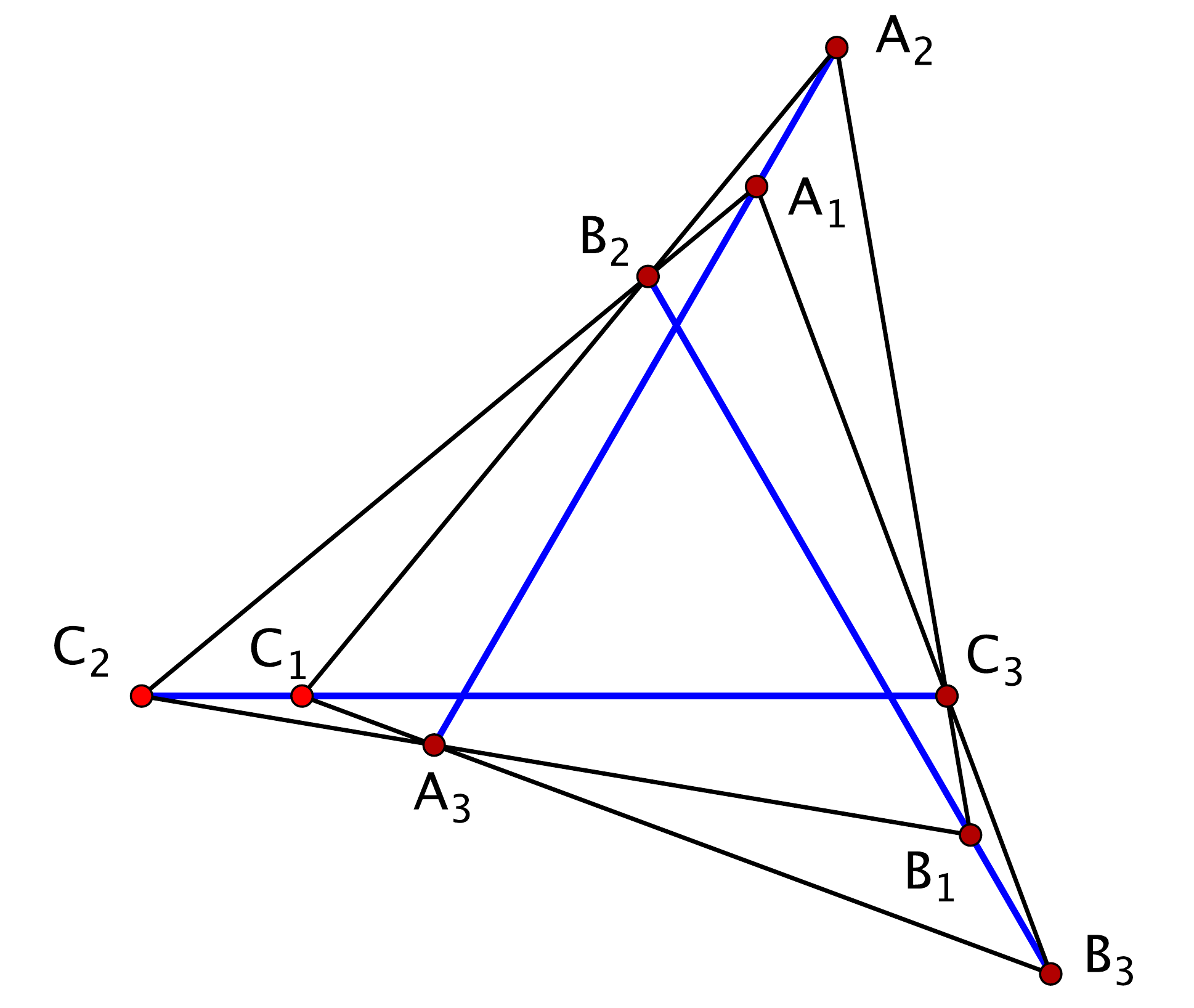}
\caption{A symmetric realization of the marked boxes $i(\Theta), \tau_1(\Theta), \tau_2(\Theta)$.}
\label{threefold}
\end{figure}

\proof
For the proof of the first statement, one can realize the box $\Theta$ in such a way that the three-fold rotational symmetry is manifestly present, see Figure \ref{threefold}. Namely,
\begin{equation*}
\begin{split}
&\Theta = (B_3,B_1,A_1,A_3;B_2,A_2),\ i(\Theta) = (A_1,A_3,B_1,B_3;A_2,B_2),\\
&\tau_1(\Theta) = (B_1,B_3,C_1,C_3;B_2,C_2),\ \tau_2(\Theta) = (C_1,C_3,A_1,A_3;C_2,A_2).
\end{split}
\end{equation*}

In terms of the marked box coordinates $(x,y)$, described in (\ref{invol}), the three operations, $i, \tau_1$, and $\tau_2$, act in the same way: $[x,y] \mapsto [1-y,x].$

\begin{figure}[hbtp]
\centering
\includegraphics[height=2in]{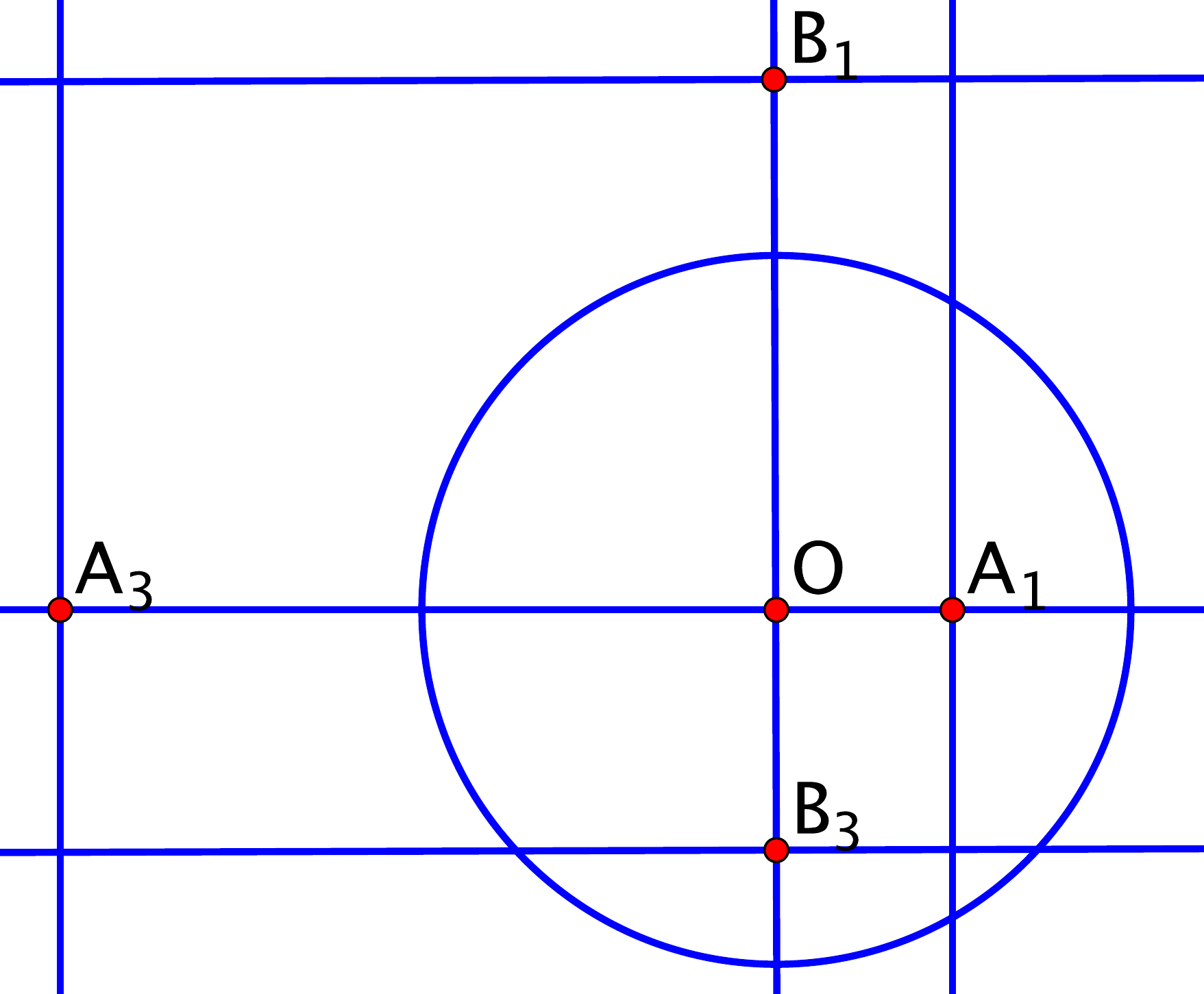}
\caption{Projective equivalence of $i(\Theta)$ and $\Theta^*$.}
\label{inversion}
\end{figure}

For the second statement, consider another realization depicted in Figure \ref{inversion}. The marked points $A_2$ and $B_2$ are at infinity, and $|OA_1| |OA_3| = |OB_1| |OB_3| = 1$. Then the polarity with respect to the unit circle centered at point $O$ acts as follows:
$$
A_1 \mapsto A_3B_2,\ A_3 \mapsto A_1 B_2,\ B_1 \mapsto B_3 A_2,\ B_3 \mapsto B_1 A_2,
$$
providing the desired projective equivalence.
\proofend

If one identifies the projective plane with its dual by a polarity, then the above discussion describes
a faithful representation of the modular group $M$ as the group of projective symmetries of 
the $G$-orbit $\Omega$ of a convex marked box. 

A marked box $\Theta$ determines a natural map $f$ of the set of vertices of the graph $\Gamma$ to the set of the marked points of the orbit $\Omega$. The map $f$ conjugates the actions of the group $G'$ on the graph $\Gamma$ and the group $M'$ of projective symmetries of the orbit $\Omega$.
The set of vertices of $\Gamma$ is dense on the circle at infinity of the hyperbolic plane $S^1$, see Figure \ref{Farey}. Using the nested properties of the interiors of the boxes in $\Omega$ and estimates on their sizes (in the elliptic plane metric), Schwartz proves the following theorem.

\begin{theorem} \label{curve}
The map $f$ extends to a homeomorphism of $S^1$ to its image. 
\end{theorem}
 
The image $\Lambda = f(S^1)$ is called the {\it Pappus curve}; see Figure \ref{Pappusiter} that provides an approximation of this curve.  

The above discussion shows that the Pappus curve is projectively self-similar. In the exceptional case of $x=y=1/2$, the curve $\Lambda$ is a straight line. Otherwise, it is not an algebraic curve, see \cite{Ha}.

The tops and bottoms of the marked boxes form a countable collection of lines that also extends to a continuous family, a curve $L$ in the dual projective plane. 

Define a transverse line field along $\Lambda$ as a continuous family of lines such that each line from the family intersects the curve at exactly one point and every point of $\Lambda$ is contained in some line. 
\begin{theorem} \label{transverse}
If the Pappus curve $\Lambda$ is not a straight line, then $L$ is a unique transverse line field along $\Lambda$.
\end{theorem}

This theorem, the fact that the Pappus curve is projectively self-similar, and computer experiments suggest that $\Lambda$ is a true fractal (unless it is a straight line). The thesis \cite{Ki} contains some preliminary numerical results on the box dimension of the Pappus curve and its dependence on the coordinates $[x,y]$ of the initial convex marked box. According to these experiments, the maximal possible box dimension of $\Lambda$ is about 1.25.

Finding the fractal dimensions of the Pappus curves as a function of $[x,y]$ or, at least, proving that this dimension is greater than one in all non-exceptional cases $[x,y]\neq [1/2,1/2]$, is an outstanding open problem. 

\section{Steiner theorem and the twisted cubic} \label{itSteiner}

This section is based on another recent ramification of the Pappus theorem, the work of J. F. Rigby \cite{Ri} and P. Hooper \cite{Ho}. 

\begin{figure}[hbtp]
\centering
\includegraphics[height=2in]{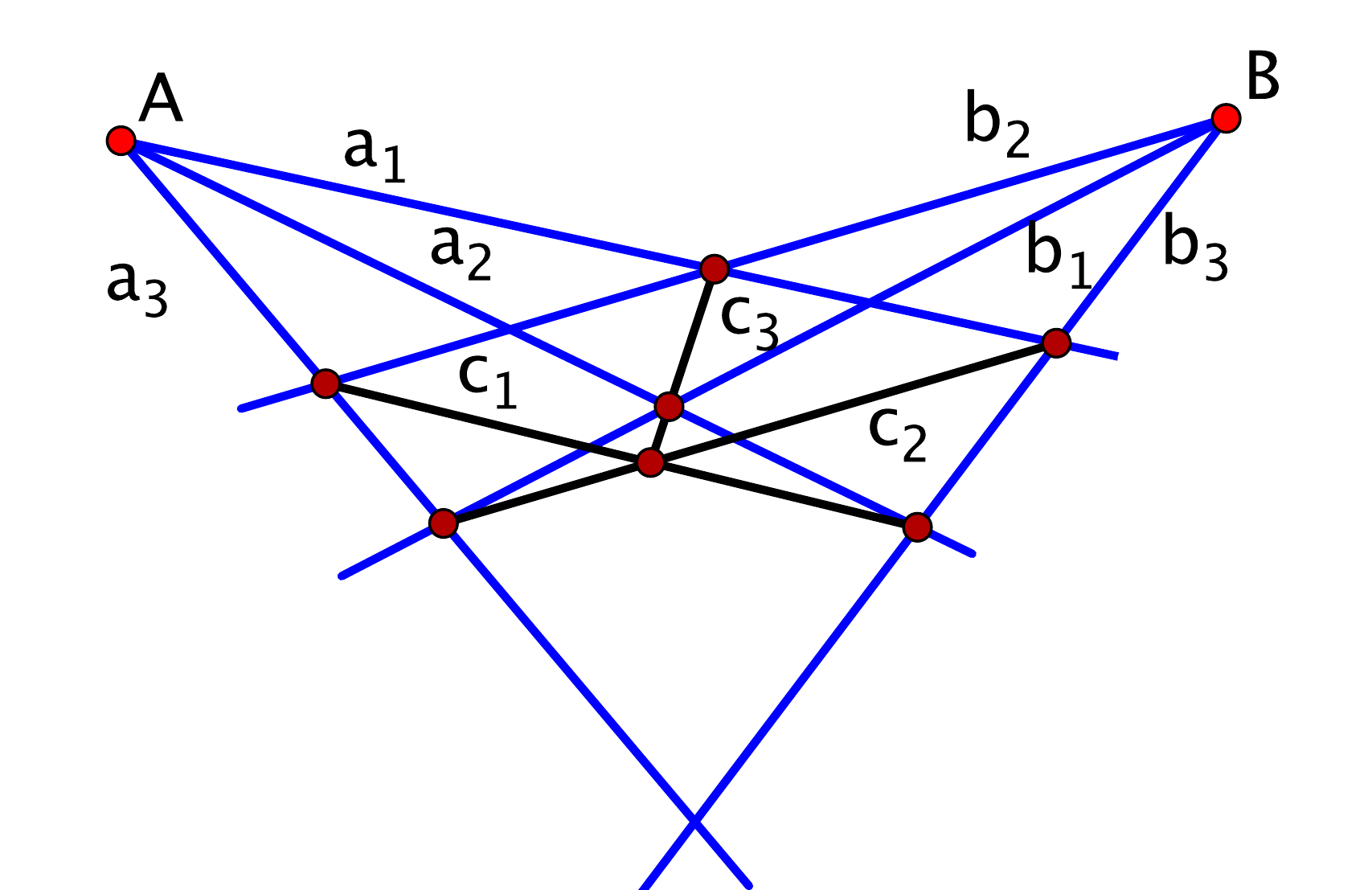}
\caption{Dual Pappus theorem.}
\label{dualPap}
\end{figure}

Let us start with the dual Pappus theorem, see Figure \ref{dualPap} where the objects dual to the ones in Figure \ref{Pappusfig} are denoted by the same letters, with the upper and lower cases swapped (the Pappus theorem is equivalent to its dual). As an aside, let us mention that the dual Pappus theorem has an interpretation in the theory of webs: the 3-web, made of three pencils of lines, is flat, see \cite{FT}, lecture 18.

Now consider Pascal's theorem, Figure \ref{Pascalfig}. The six permutations of the points on the conic yield 60 Pascal lines. These lines and their intersection points, connected by further lines, form a intricate configuration of 95 points and 95 lines, the {\it hexagrammum mysticum}. There is a number of theorems describing this configuration, due to Steiner, Pl\"ucker, Kirkman, Cayley, and Salmon. See \cite{CR1,CR2} for a contemporary account of this subject.

\begin{figure}[hbtp]
\centering
\includegraphics[height=3.2in]{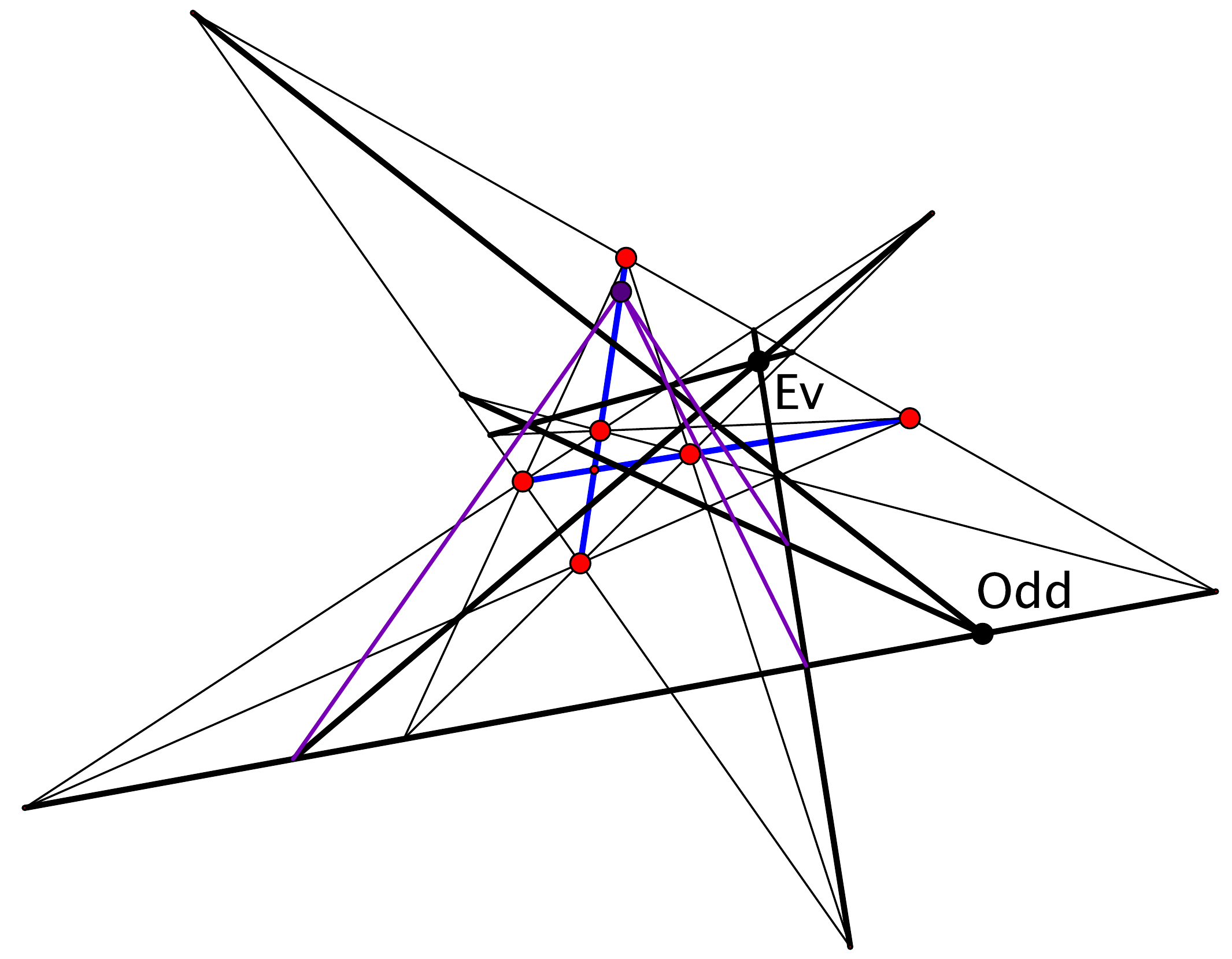}
\caption{Two Steiner points, corresponding to even and odd permutations, are labelled. One of the points $\ell^*(\varphi, s(\psi))$ is shown.}
\label{return2}
\end{figure}

The Pappus theorem is a particular case of Pascal's theorem, and in this case, the number of lines that result from permuting the initial points (say, points $B_1,B_2,B_3$ in Figure \ref{Pappusfig}) reduces to six,  shown in  Figure \ref{return2}.


Let us introduce notations. Consider Figure \ref{Pappusfig} and denote the triples of points:
$$
{\cal A} = (A_1,A_2,A_3),\ {\cal B} = (B_1,B_2,B_3).
$$
The Pappus theorem produces a new triple, ${\cal C} = (C_1,C_2,C_3)$.
The lines containing these triples are denoted by $a,b,c$, respectively.
We write: $c = \ell ({\cal A},{\cal B}).$ 

We use a similar notation for the dual Pappus theorem: if 
$$
\alpha=(a_1,a_2,a_3),\ \beta=(b_1,b_2,b_3) 
$$
are two triples of concurrent lines, then $\ell^* (\alpha,\beta)$ is the point of intersection of the triple of lines $(c_1,c_2,c_3)$, see Figure \ref{dualPap}.

The permutation group $S_3$ acts on triples by the formula
$$
s({\cal B}) = (B_{s^{-1}(1)},B_{s^{-1}(2)},B_{s^{-1}(3)}).
$$
Let $\sigma\in S_3$ be a cyclic permutation, and $\tau\in S_3$ be a transposition of two elements.

The following result, depicted in Figure \ref{return2}, is due to J. Steiner.

\begin{theorem} \label{St1}
The three Pappus lines $\ell ({\cal A},s({\cal B}))$ where $s\in S_3$ is an even permutation, are concurrent, and so are the three lines corresponding to the odd permutations.
\end{theorem} 

Thus we obtain two triples of concurrent lines; denote them by
$$
\varphi = (\ell ({\cal A},{\cal B}), \ell ({\cal A},\sigma({\cal B}), \ell ({\cal A},\sigma^2({\cal B})), 
\psi = (\ell ({\cal A},\tau({\cal B})), \ell ({\cal A},\tau\sigma({\cal B}), \ell ({\cal A},\tau\sigma^2({\cal B})).
$$ 

Apply the dual Pappus theorem to  the permutations of these triples of lines. By the dual Steiner theorem, the six points $\ell^*(\varphi, s(\psi)),\ s\in S_3$, are collinear in threes. 

More is true. The next two theorems are due to Rigby \cite{Ri}.

\begin{theorem} \label{Ri1}
The points $\ell^*(\varphi, s(\psi))$ lie on line $a$ when $s$ is an even permutation, and on line $b$ when $s$ is odd.
\end{theorem} 

Let ${\cal B}'$ be another collinear triple of points such that the line $b'$ still passes through point $O = a\cap b$. Applying the above constructions to ${\cal A}, {\cal B}'$, we obtain new triples of lines $\varphi', \psi'$, and a new triple of points $\ell^*(\varphi', s(\psi'))$ on line $a$ where $s$ is an even permutation.

\begin{theorem} \label{Ri2}
The new triple of points coincides with the old one: for even permutations $s$, one has $\ell^*(\varphi', s(\psi'))= \ell^*(\varphi, s(\psi))$.
\end{theorem}

Theorems \ref{Ri1} and \ref{Ri2} are stated by Rigby without proof; to quote,
\begin{quote}
The theorems in this section have been verified in a long and tedious manner using coordinates. There seems little point in publishing the calculations; it is to be hoped that shorter and more elegant proofs will be found.
\end{quote}
Conceptual proofs are given in \cite{Ho}; the reader is referred to this paper and is encouraged to find an alternative approach to these results.   

The above theorems make it possible to define the {\it Steiner map}
$$
S_O: (A_1,A_2,A_3) \mapsto (\ell^*(\varphi, \psi), \ell^*(\varphi, \sigma^2(\psi)), \ell^*(\varphi, \sigma(\psi))). 
$$ 
This map depends on the point $O$, but not on the choice of the triple ${\cal B}$. 

The Steiner map commutes with permutations of the points involved, and hence it induces a map of the space of unordered triples of points of the projective line. Abusing notation, we denote this induced map by the same symbol.
Hooper \cite{Ho} gives a complete description of the Steiner map. 

Assume that the ground field is  complex. The space of unordered triples of points of $\CP^1$, that is, the symmetric cube $S^3(\CP^1)$, is identified with $\CP^3$. This is a particular case of the Fundamental Theorem of Algebra, one of whose formulations is that $S^n(\CP^1)=\CP^n$ (given by projectivizing the Vieta formulas that relate the coefficients of a polynomial to its roots). Thus $S_O$ is a self-map of $\CP^3$.

The set of cubic polynomials with a triple root corresponds to a curve $\Gamma \subset \CP^3$, the twisted cubic (the moment curve). The secant variety of the twisted cubic, that is, the union of its tangent and secant lines, covers $\CP^3$, and the lines are pairwise disjoint, except at the points of $\Gamma$. 

The set of cubic polynomials with a zero root corresponds to a plane in $\CP^3$. Denote this plane by $\Pi$. The Steiner map $S_O: \CP^3 \to \CP^3$ is described in the next theorem.

\begin{theorem} \label{Stmap}
(i) The map $S_O$ preserves the secants of the twisted cubic $\Gamma$ that do not pass through the origin (the image of the cubic polynomial $z^3$). \\
(ii) One can choose projective coordinates on these secant lines so that the map is given by the formula $x \mapsto x^2$. \\
(iii) The choice of coordinates is as follows:  the two points of intersection of the secant line with $\Gamma$ have coordinates $0$ and $\infty$, and the intersection point of the secant with the plane $\Pi$ has coordinate $-1$.
\end{theorem}

In homogeneous coordinates of $\CP^3$, the map $S_O$ is polynomial of degree 6; see \cite{Ho} for an explicit formula for a particular choice $O=(0:1)$. 

In the real case, the secant lines are identified with the circle $\R/\Z$ , and the Steiner map becomes the doubling map $t \mapsto 2t$ mod 1, a well known measure preserving ergodic transformation. 

\section{Pentagram-like maps on inscribed polygons} \label{pentalike}

This section, based on \cite{ST1}, concerns eight configuration theorems of projective geometry that were discovered in the study of the pentagram map. 

\begin{figure}[hbtp]
\centering
\includegraphics[height=1.7in]{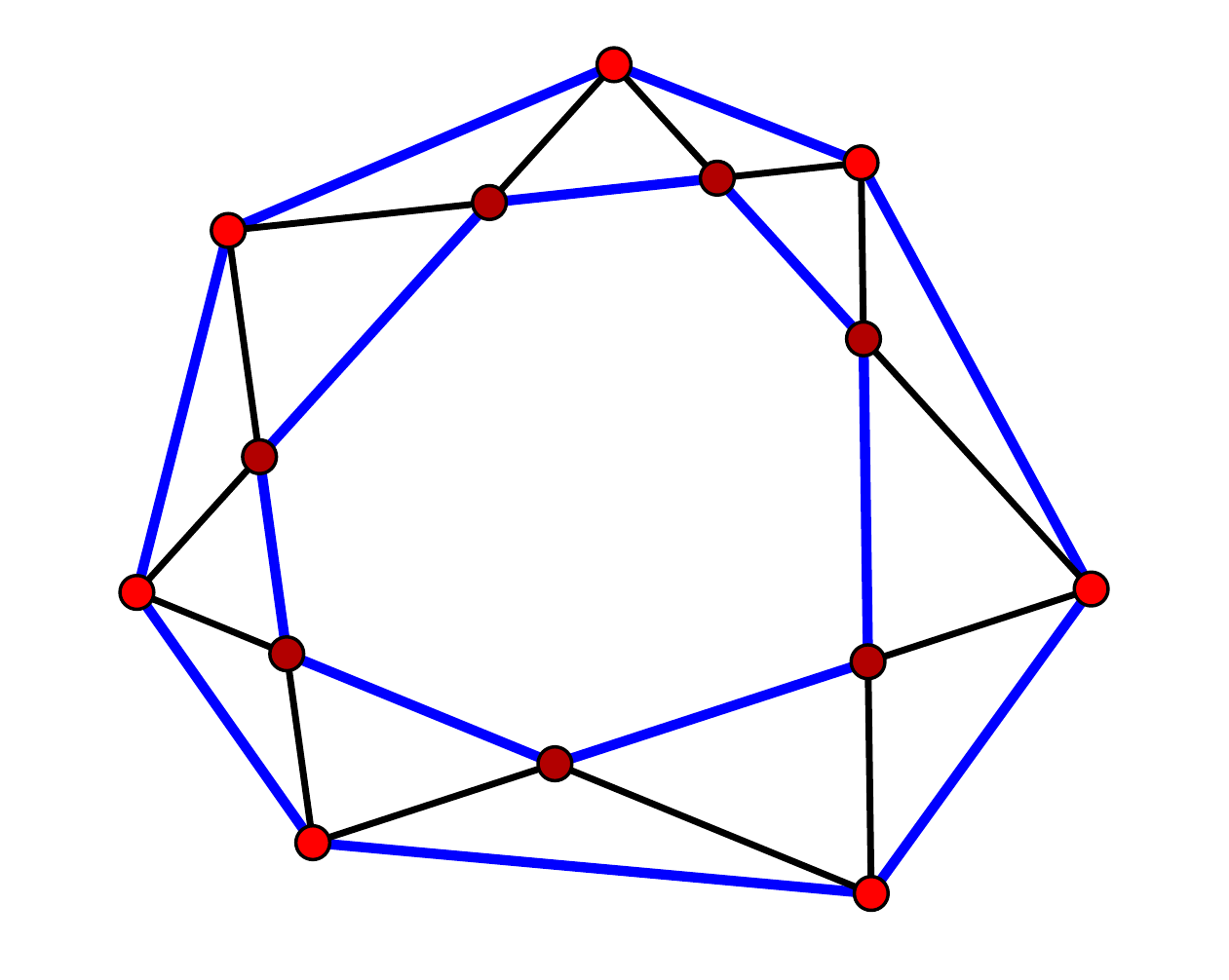}
\caption{The pentagram map takes an $n$-gon $P$ to the polygon made by the intersection points of the short (skip one) diagonals of $P$. Here $n=7$.}
\label{penta}
\end{figure}

The pentagram map, whose study was put forward by R. Schwartz  \cite{Sch92}, is a transformation of the moduli space of projective equivalence classes of polygons in the projective plane depicted in Figure \ref{penta}. The pentagram map has become a popular object of study: it is a discrete completely integrable system, closely related with the  theory of cluster algebras. See \cite{GSTV,Gl,GP,OST1,OST2,So} for a sampler of the current literature on this subject.

To formulate the results, let us introduce some notations.

By a polygon in the projective plane we mean a cyclically ordered collection of its vertices (that also determines the cyclically ordered collection of lines, the sides of the polygon).

Let ${\cal C}_n$ and ${\cal C}_n^*$ be the spaces $n$-gons in the projective plane $\RP^2$ and its dual $(\RP^2)^*$. Define the $k$-diagonal map 
$T_k: {\cal C}_n  \to  {\cal C}_n^*$: for $P=\{p_1,...,p_n\}$, 
$$T_k(P)=\{(p_1p_{k+1}),(p_2p_{k+2}),\ldots, (p_np_{k+n})\}.$$
Each map  $T_k$ is an involution; the map $T_1$ is the projective duality that sends a polygon to the cyclically ordered collection of its sides.

Extend the notation to muti-indices: $T_{ab}=T_a \circ T_b, T_{abc}=T_a \circ T_b\circ T_c$, etc.
For example, the  pentagram map is $T_{12}$. If $P$ is a polygon in $\RP^2$ and $Q$ a polygon in $(\RP^2)^*$, and there exists a projective transformation $\RP^2 \to (\RP^2)^*$ that takes $P$ to $Q$, we write: $P \sim Q$.

Now we are ready to formulate our results; they  concern polygons inscribed into a conic or circumscribed about a conic. 

\begin{theorem} \label{thm12}
(i) If $P$ is an inscribed $6$-gon, then $P \sim T_2(P)$.\\
(ii) If $P$ is an inscribed $7$-gon, then $P \sim T_{212}(P)$.\\
(iii) If $P$ is an inscribed $8$-gon, then $P \sim T_{21212}(P)$.
\end{theorem}

\begin{figure}[hbtp]
\centering
\includegraphics[height=2in]{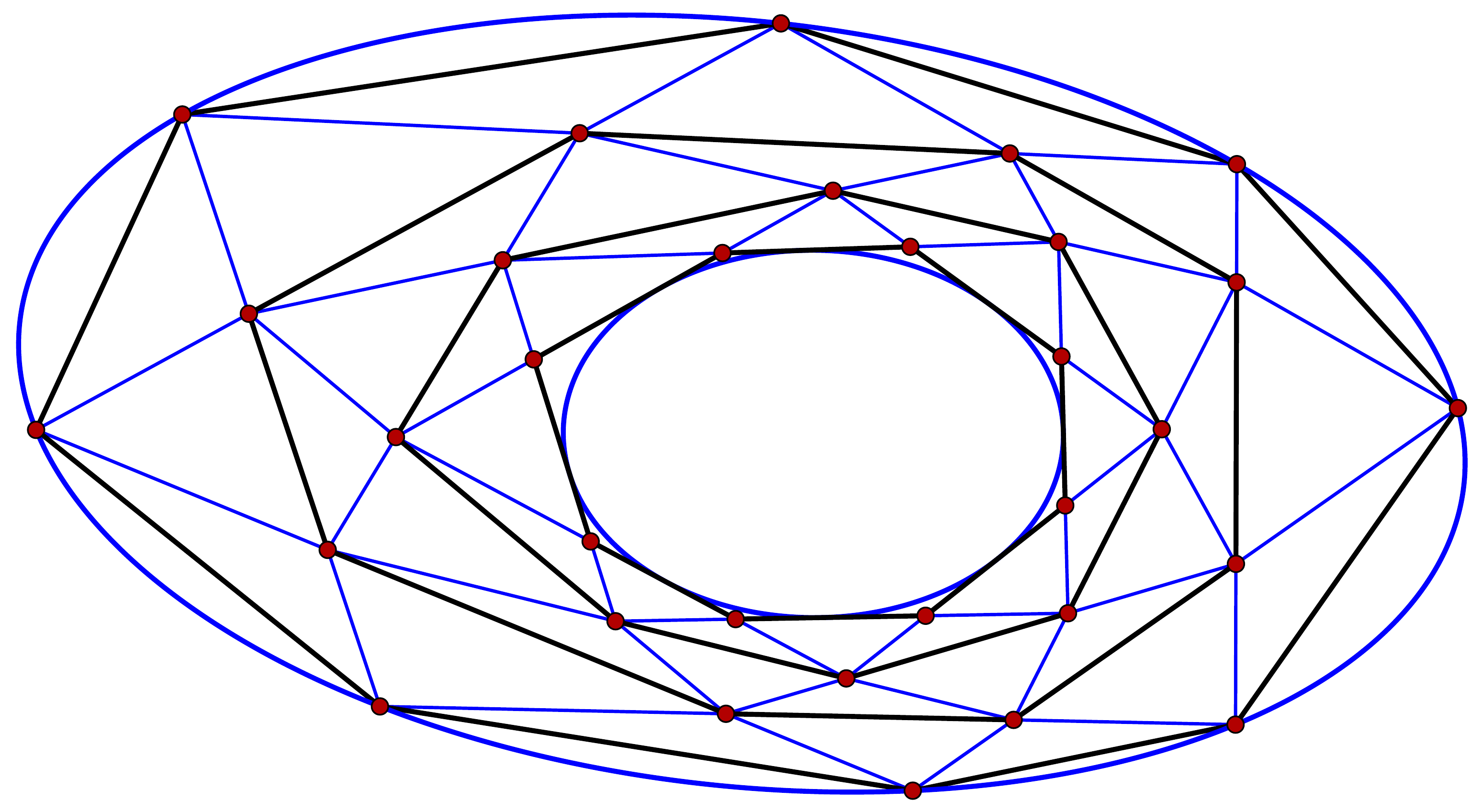}
\caption{The third iteration of the pentagram map on an inscribed octagon yields a projectively dual octagon.}
\label{octagon}
\end{figure}

Surprisingly, this sequence doesn't continue!
Theorem \ref{thm12} (iii) is depicted in Figure \ref{octagon}. See also Schwartz's applet \cite{App1} for illustration of this and other results of this section.

\begin{theorem} \label{thm313}
If $P$ is a circumscribed $9$-gon, then $P \sim T_{313}(P)$.
\end{theorem}

\begin{figure}[hbtp]
\centering
\includegraphics[height=2.3in]{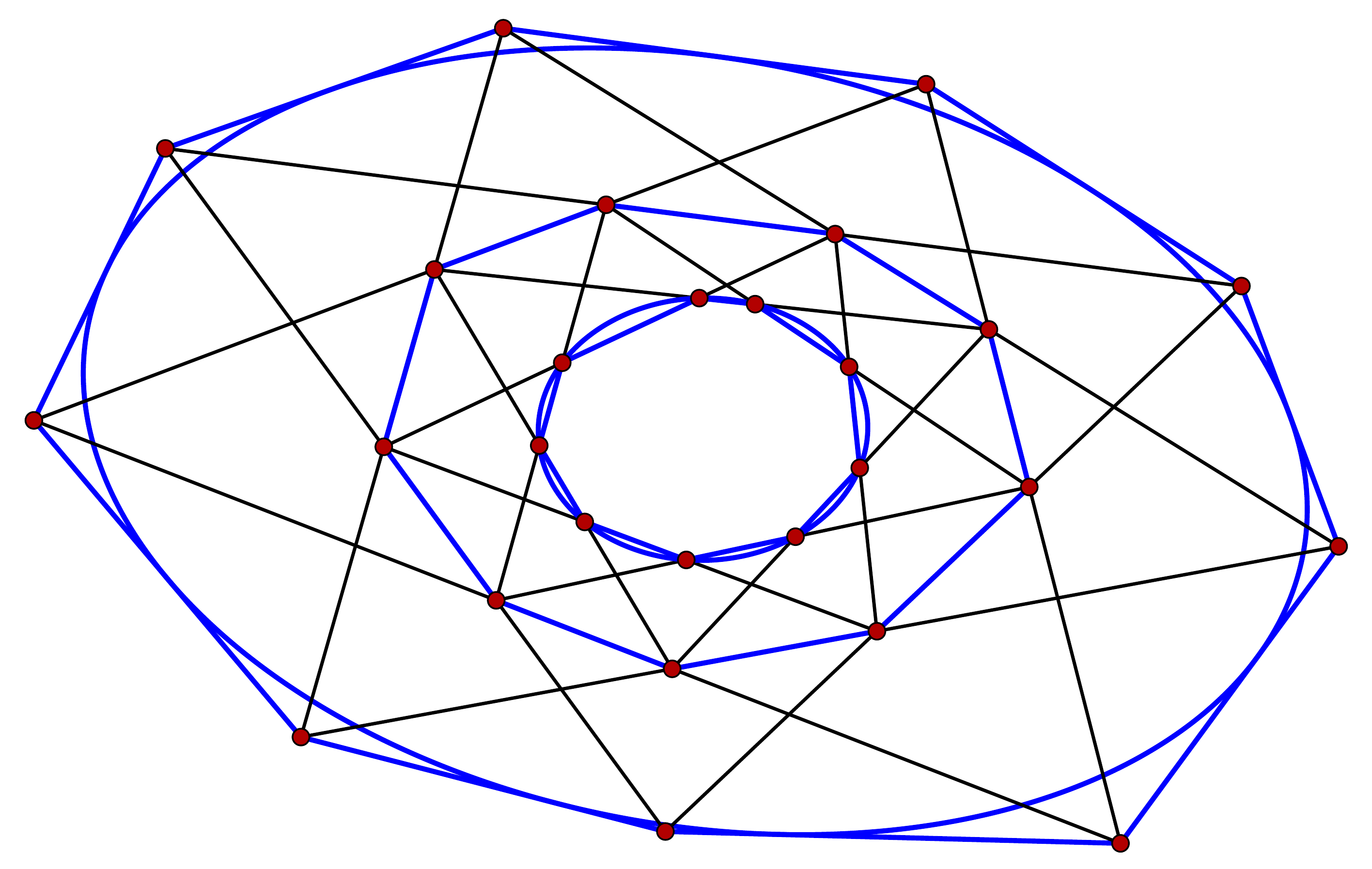}
\caption{Theorem \ref{thm313}.}
\label{nonagon}
\end{figure}

See Figure \ref{nonagon}.

\begin{theorem} \label{thm34}
If $P$ is an inscribed $12$-gon, then $P \sim T_{3434343}(P)$.
\end{theorem}

The next results have a somewhat different flavor: one does not claim anymore that the final polygon is projectively related to the initial one.

\begin{theorem} \label{thm13}
(i) If $P$ is an inscribed $8$-gon, then $T_{3}(P)$ is circumscribed.\\
(ii) If $P$ is an inscribed $10$-gon, then $T_{313}(P)$ is circumscribed.\\
(iii)  If $P$ is an inscribed $12$-gon, then $T_{31313}(P)$ is circumscribed.
\end{theorem}

\begin{figure}[hbtp]
\centering
\includegraphics[height=2.5in]{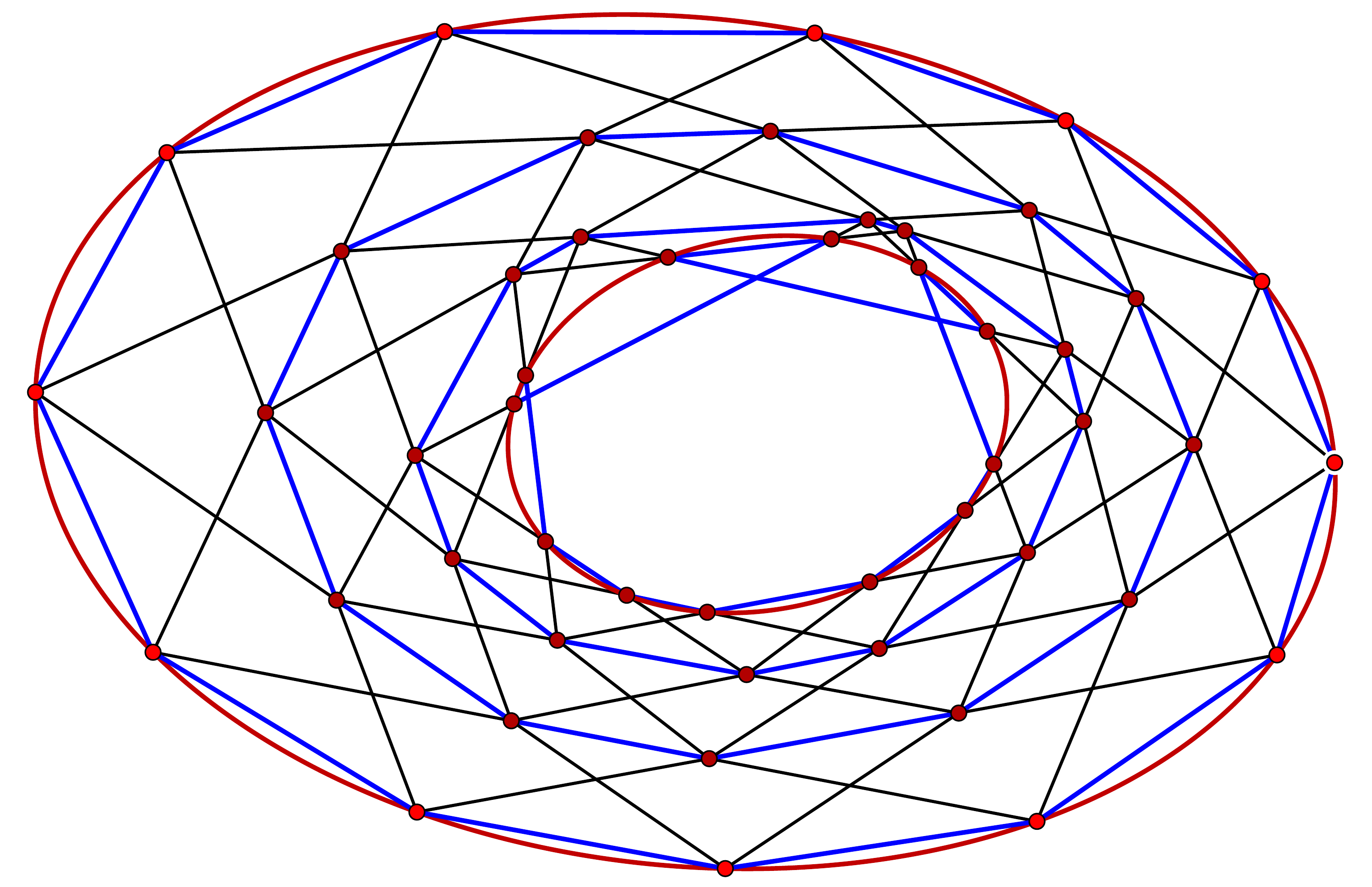}
\caption{Theorem \ref{thm13} (iii).}
\label{dodecagon}
\end{figure}

Again, in spite of one's  expectation, this sequence does not continue. Theorem \ref{thm13} (iii) is illustrated in Figure \ref{dodecagon}.

Now about the discovery of these results and their proofs. Theorems \ref{thm12} (i) and (ii) were discovered in our study of the pentagram map. Then V. Zakharevich, a participant of the 2009 Penn State REU (Research Experience for Undergraduates) program, discovered Theorem \ref{thm313}. Inspired by this discovery, we did an extensive computer search for this kind of configuration theorems; the results are the above eight theorems. We think that the list above is exhaustive, but this remains a conjecture.

Note that one may cyclically relabel the vertices of a polygon to deduce seemingly new  theorems. Let us illustrate this by an example. Rephrase the statement of Theorem \ref{thm13} (iii) as follows: {\it If $P$ is an inscribed dodecagon then $T_{131313}(P)$ is also inscribed}. Now relabel the vertices by $\sigma(i)=5i$ mod 12. The map $T_3$ is conjugated by $\sigma$ as follows: 
$$
i\mapsto 5i\mapsto 5i+3\mapsto 5(5i+3)=i+3\ \ {\rm mod}\ 12, 
$$
that is, it is the map is $T_3$ again, and the map $T_1$ becomes
$$
i\mapsto 5i\mapsto 5i+1\mapsto 5(5i+1)=i+5\ \ {\rm mod}\ 12,
$$
that is, the map is $T_5$. One arrives at the statement: {\it If $P$ is an inscribed dodecagon then $T_{535353}(P)$ is also inscribed}. 

We proved all of the above theorems, except Theorem \ref{thm13} (iii), by uninspiring computer calculations (the  symbolic manipulation required for a proof of Theorem \ref{thm13} (iii) was beyond what we could manage in Mathematica). 

Of course, one wishes for elegant geometric proofs. Stephen Wang found proofs of Theorems \ref{thm13} (i) and (ii) which are presented below, and Maria Nastasescu, a 2010 Penn State REU participant, found algebraic geometry proofs of the same two theorems. Fedor Nilov proved Theorem \ref{thm13} (iii) using a planar projection of hyperboloid of one sheet. Unfortunately, none of these proofs were published.

Here is Wang's proof of  Theorem \ref{thm13} (i). 

\begin{figure}[hbtp]
\centering
\includegraphics[height=2.2in]{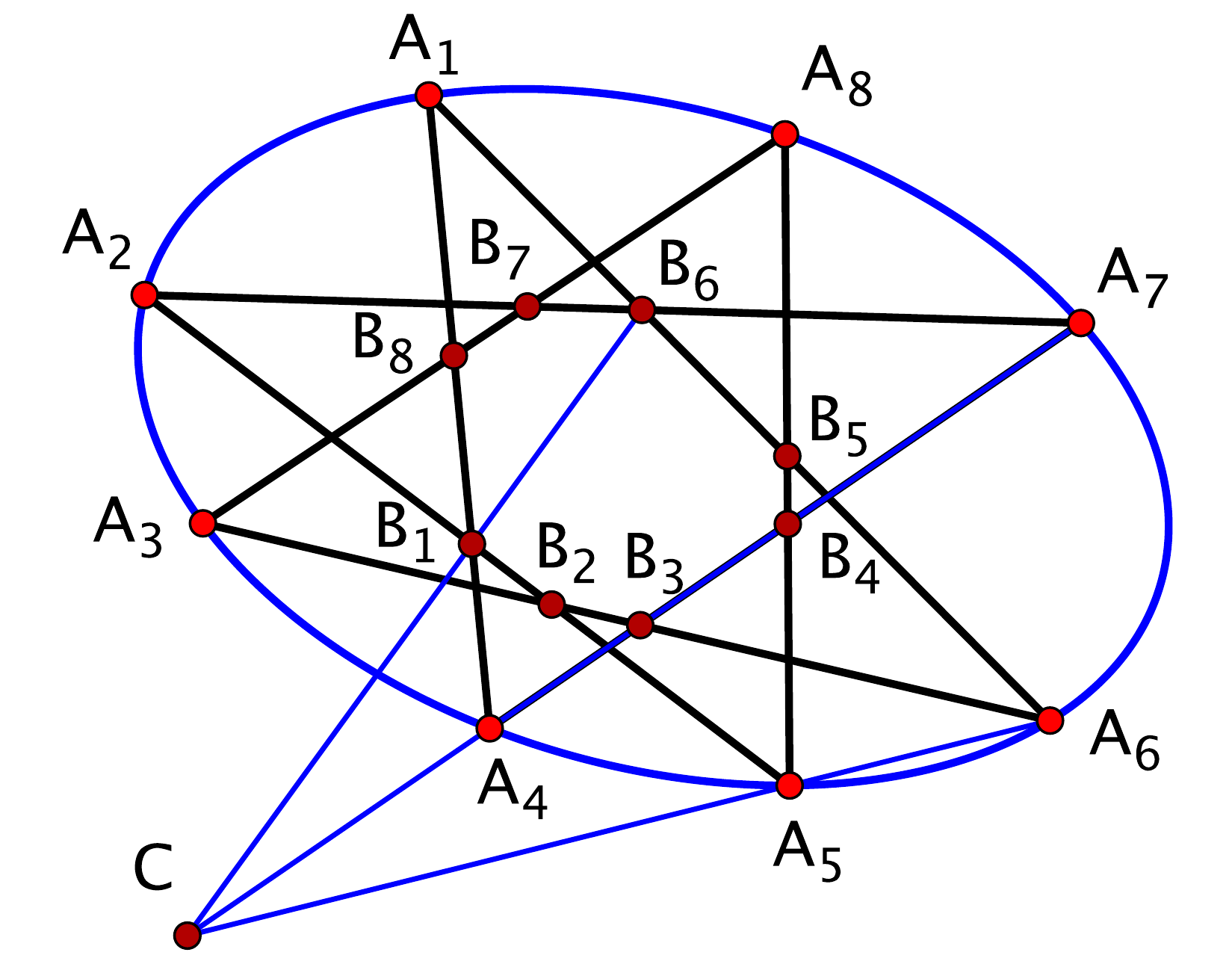}
\caption{Proof of Theorem \ref{thm13} (i).}
\label{Wang2}
\end{figure}

Consider Figure \ref{Wang2}. We need to prove that the points $B_1,\ldots, B_8$ lie on a conic.

The hexagon $A_6 A_1 A_4 A_7 A_2 A_5$ is inscribed so, by Pascal's theorem, the points $B_1, B_6$ and $C$ are collinear. That is, the intersection points of the opposite sides of the hexagon $B_1 B_2 B_3 B_4 B_5 B_6$ are collinear. By the converse Pascal theorem, this hexagon is inscribed. 

A similar argument shows that the hexagon $B_2 B_3 B_4 B_5 B_6 B_7$ is inscribed. But the two hexagons share five vertices, hence they are inscribed in the same conic. Likewise, $B_8$ lies on this conic as well.

\begin{figure}[hbtp]
\centering
\includegraphics[height=3.4in]{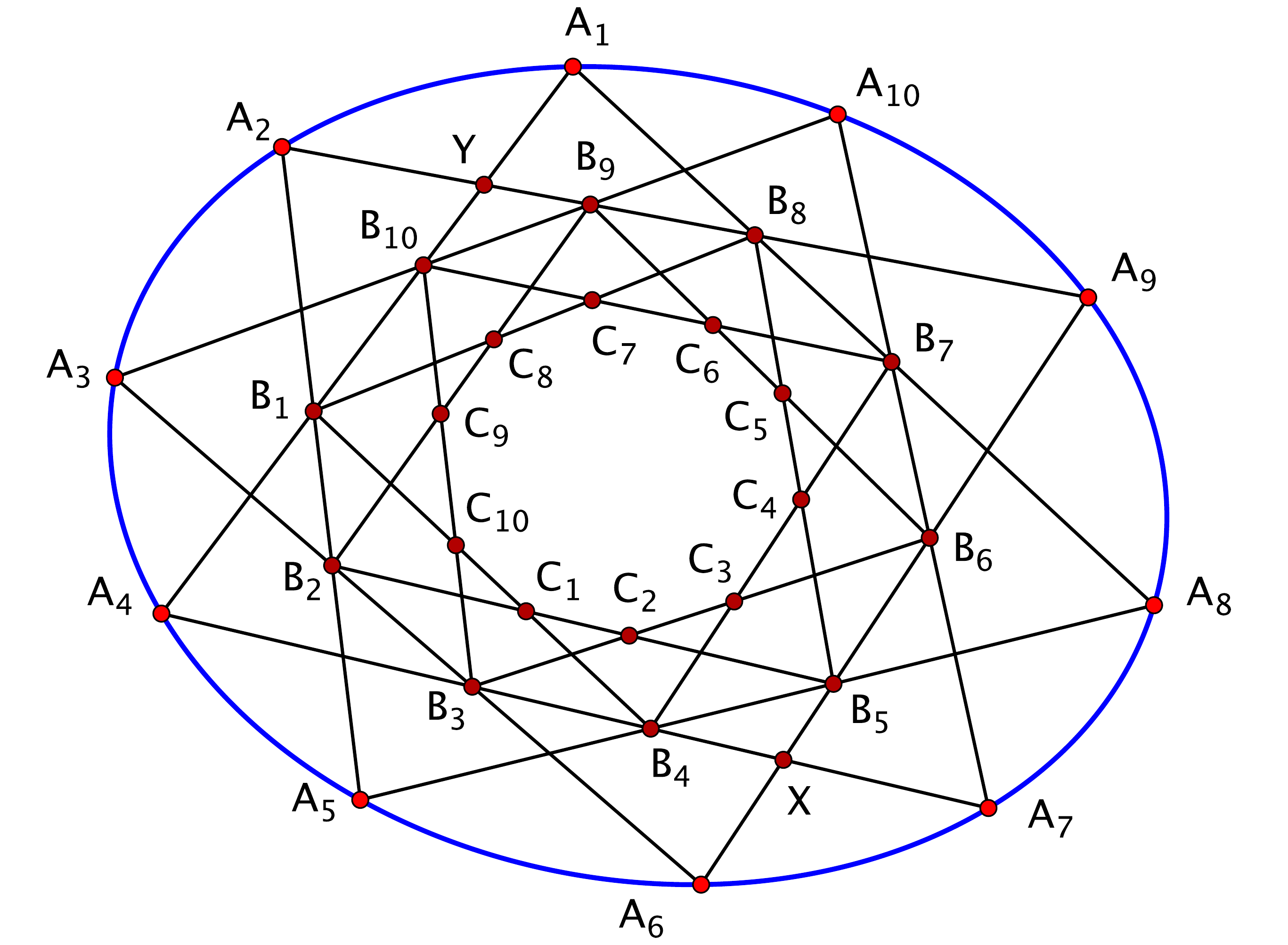}
\caption{Proof of Theorem \ref{thm13} (ii).}
\label{Wang3}
\end{figure}

Now, to the proof of Theorem \ref{thm13} (ii), see Figure \ref{Wang3}.

Consider the inscribed hexagon $A_3A_6A_9A_{10}A_7A_4$. By Pascal's theorem, the points
$$
(A_3A_6)\cap (A_{10}A_7),\ (A_6A_9)\cap (A_7A_4),\ (A_9A_{10})\cap (A_4A_3)
$$
are collinear. Hence the triangles $A_3B_3A_4$ and $A_{10}A_9B_6$ 
are perspective. By the Desargues theorem, the points
$ A_4, A_9, (B_3B_6)\cap (A_3A_{10})$ are collinear. 

It follows that the triangles $B_9B_{10}Y$ and $B_3B_6X$ are perspective. 
By the Desargues theorem, the points $X,Y, (B_6B_9)\cap (B_3B_{10})$ are collinear. 
The same argument, with all indices shifted by five, implies that the points 
$X,Y, (B_1B_4)\cap (B_8B_5)$ are collinear as well. Hence the points 
$$
(B_6B_9)\cap (B_3B_{10}), (B_1B_4)\cap (B_8B_5),\ {\rm and}\ X
$$
are collinear. Reinterpret this as the collinearity of
$$
(C_{10}B_{10})\cap (C_5B_9),\ (C_{10}B_4)\cap (C_5B_5),\ (B_3B_4)\cap (B_5B_6).
$$
It follows that the triangles $B_3B_4C_{10}$ and $B_5B_6C_5$ are perspective. 
By the Desargues theorem, the points $B_4,B_5$ and $(C_{10}C_5)\cap (C_2C_3)$ are collinear. That is, the points
$$
(C_{10}C_5)\cap (C_2C_3), (C_{10}C_1)\cap (C_3C_4),\ (C_1C_2)\cap (C_4C_5)
$$
are collinear, and by the converse Pascal theorem, the points 
$$C_{10},C_1,C_2,C_3,C_4,C_5$$
lie on a conic. The rest is the same as in the previous proof. 
\proofend

One can add to Theorems \ref{thm12} -- \ref{thm13} a statement about pentagons.  Consider the following facts:\\
(i) every pentagon is inscribed in a conic and circumscribed about a conic;\\
(ii) every pentagon is projectively equivalent to its dual;\\
(iii) the pentagram map sends every pentagon to a projectively-equivalent one.\\
Therefore one may add the following theorem to our list: {\it for a pentagon $P$, one has $P\sim T_2(P)$}.

The following result of R. Schwartz \cite{Sch01,Sch08,Gl} also has a similar flavor.

\begin{theorem} \label{degen}
If $P$ is a $4n$-gon inscribed into a degenerate conic (that is, a pair of lines) then 
$$(T_1T_2T_1T_2\dots T_1)(P)\qquad (4n-3\  {\rm terms)}$$
 is also inscribed into a degenerate conic.
 \end{theorem}
 
 One wonders whether there is a unifying theme here. A possibly relevant reference is \cite{GP}.

\section{Poncelet grid, string construction, and billiards in ellipses} \label{grid}

A Poncelet polygon is a polygon that is inscribed into an ellipse $\Gamma$ and circumscribed about an ellipse $\gamma$. Let $L_1,\ldots, L_n$ be the lines containing the sides of a Poncelet $n$-gon, enumerated in such a way that their tangency points with $\gamma$ are in the cyclic order. The {\it Poncelet grid} is the collection of $n(n+1)/2$ points  $L_i \cap L_j$, where  $L_i \cap L_i$ is the tangency point of the line $L_i$ with $\gamma$. To simplify the exposition, assume that $n$ is odd (for even $n$, the formulations are slightly different).

\begin{figure}[hbtp]
\centering
\includegraphics[height=2.2in]{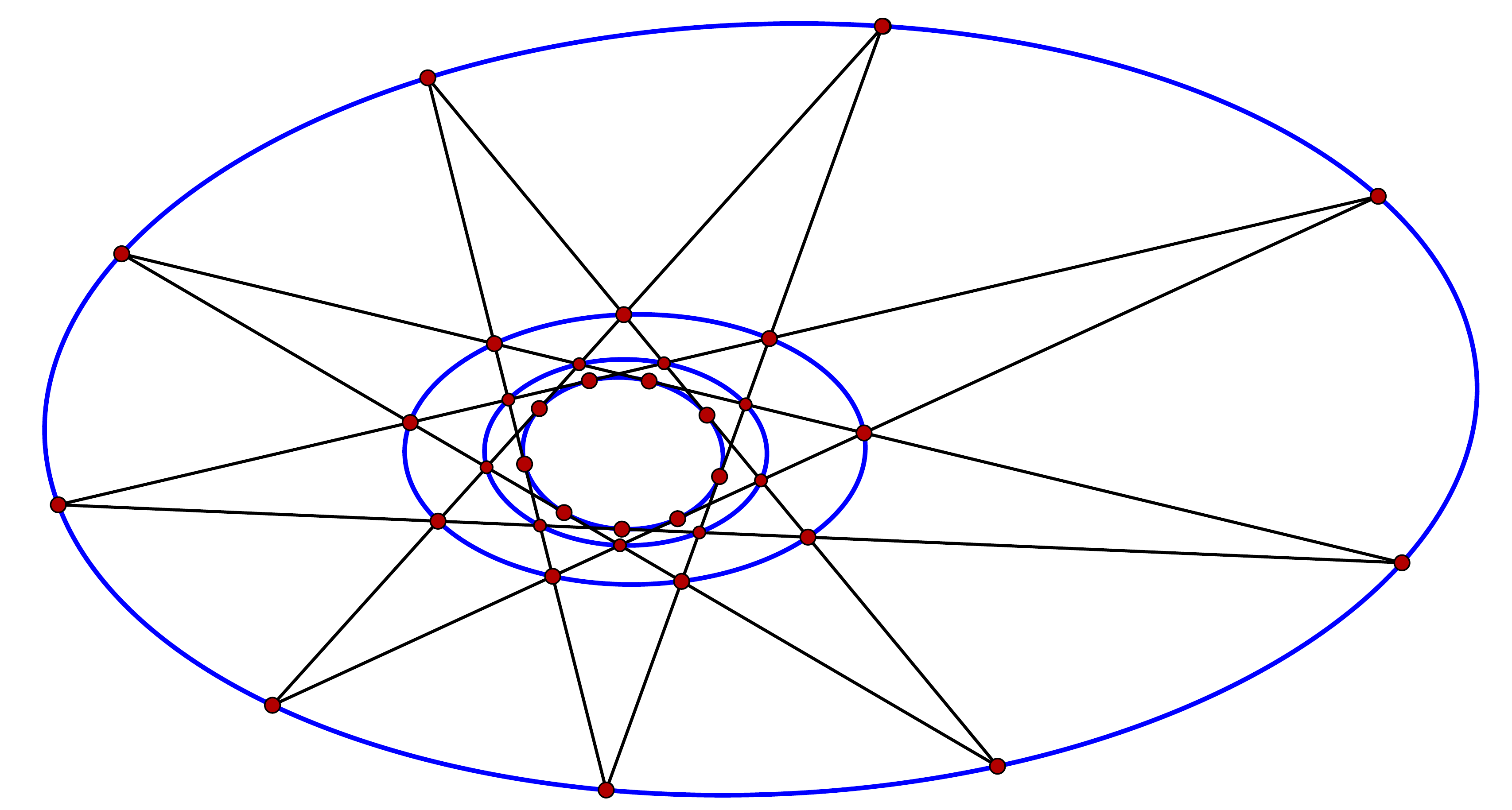}
\caption{Poncelet grid, $n=9$: shown are the concentric sets $P_0, P_2, P_3$, and $P_4$  that lie on four ellipses.}
\label{Grid}
\end{figure}

One can partition the Poncelet grid in two ways. Define the sets
$$
P_k = \cup_{i-j=k} \ell_i \cap \ell_j,\quad Q_k = \cup_{i+j=k} \ell_i \cap \ell_j,
$$
where the indices are understood mod $n$.
There are $(n + 1)/2$ sets $P_k$ , each containing $n$ points, and $n$ sets $Q_k$ , each containing $(n + 1)/2$ points.
The sets $P_k$ are called concentric, and the sets $Q_k$ are called radial, 
see Figure \ref{Grid}. 

The following theorem is proved in \cite{Sch07}.

\begin{theorem} \label{Pgrid}
(i) The concentric sets   lie on nested ellipses, and the radial sets   lie on disjoint hyperbolas. \\
(ii) The complexified versions of these conics have four common  tangent lines. \\
(iii) All the concentric sets are projectively equivalent to each other, and so are all the radial sets. 
\end{theorem}

In this section, following \cite{LT}, we  prove this projective theorem using Euclidean geometry, namely, the billiard properties of conics. As a by-product of this approach, we establish the Poncelet theorem and prove the theorem of Reye and Chasles on inscribed circles. See \cite{DR,Fl} for general information about  the Poncelet theorem, and \cite{KT,Tab95,Tab05} for the theory of billiards. 

The reduction to billiards goes as follows.
Any pair of nested ellipses $\gamma \subset \Gamma$ can be taken to a pair of confocal ellipses by a suitable projective transformation. This transformation takes a Poncelet polygon to a periodic billiard trajectory in $\Gamma$.

The billiard inside a convex domain with smooth boundary is a transformation of the space of oriented lines (rays of light) that intersect the domain: an incoming billiard trajectory hits the boundary (a mirror) and optically reflects so that the angle of incidence equals the angle of reflection.

The space of oriented lines has an area form, preserved by the optical reflections (independently of the shape of the mirror). 
Choose an origin, and introduce coordinates $(\alpha,p)$ on the space of rays: $\alpha$ is the direction of the ray, and $p$ is its signed distance to the origin. Then the invariant area form is  $\omega=d\alpha\wedge d p$.

A {\it caustic} of a billiard is a curve $\gamma$ with the  property that if a segment of a billiard trajectory is tangent to $\gamma$, then so is each reflected segment. 

There is no general method of describing caustics of a given billiard curve,\footnote{The existence of caustics for strictly convex and sufficiently smooth billiard curves is proved in the framework of the KAM theory.}  but the converse problem, to reconstruct a billiard table $\Gamma$ from its caustic $\gamma$, has a simple solution given by the following {\it string construction}: wrap a non-stretchable closed string around $\gamma$, pull it tight, and move the farthest point around $\gamma$; the trajectory of this point is the billiard curve $\Gamma$. This construction yields a 1-parameter family of billiard tables sharing the caustic $\gamma$: the parameter is the length of the string.

\begin{figure}[hbtp]
\centering
\includegraphics[height=1.7in]{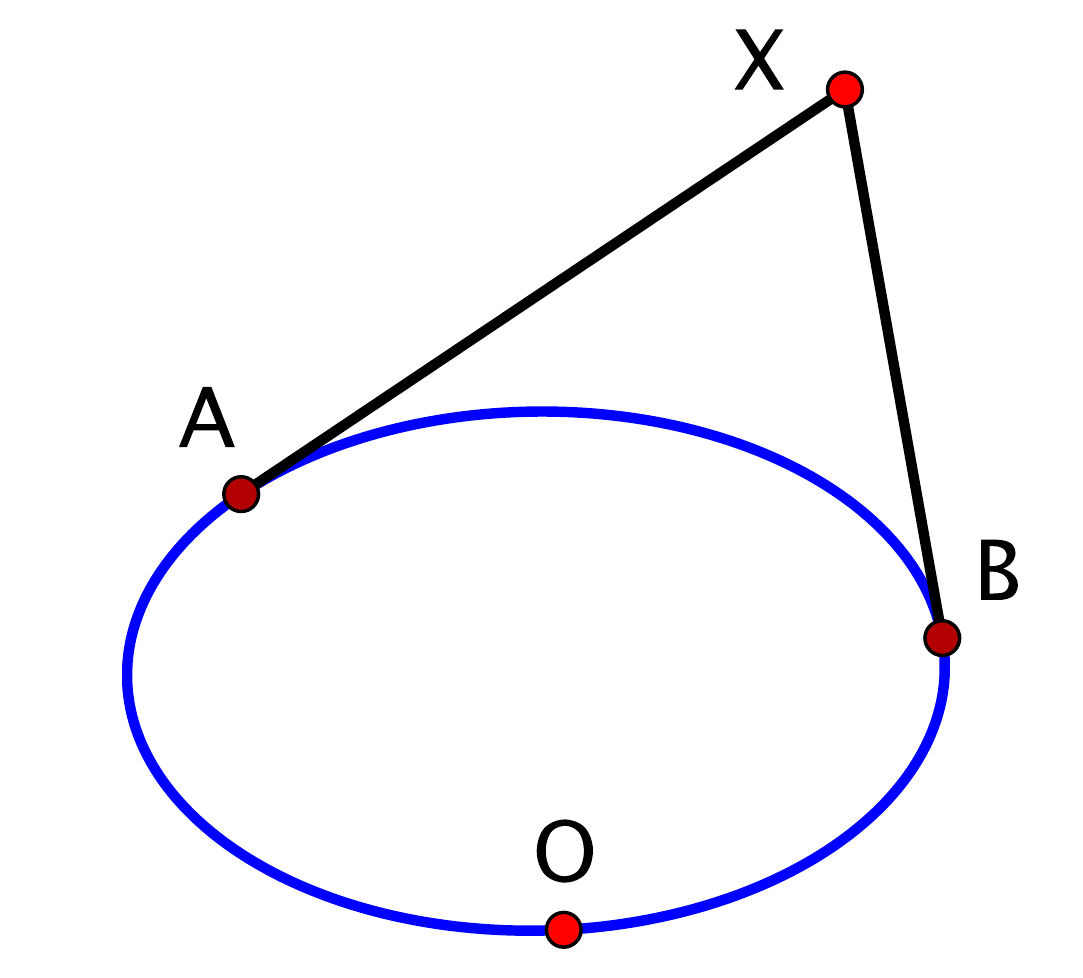}
\caption{String construction}
\label{string}
\end{figure}

The reason is as follows, see Figure \ref{string}.  For a point $X$ outside of the oval $\gamma$, consider two functions: 
$$
f(X)=|XA|+\stackrel{\smile}{|AO|},\ g(X)=|XB|+\stackrel{\smile}{|BO|}.
$$

The gradients of these functions are the unit vectors along the lines $AX$ and $BX$, respectively. It follows that these two lines make equal angles with the level curves of the functions $f+g$ and $f-g$, and that these level curves are orthogonal to each other. In particular, the level curves of $f+g$ are the billiard tables for which $\gamma$ is a caustic.

Note that the function $f+g$ does not depend on the choice of the auxiliary point $O$, whereas the function $f-g$ is  defined up to an additive constant, so its level curves are well defined.

Here is a summary of the billiard properties of conics. 
The interior of an ellipse is foliated by confocal ellipses: these are the caustics of the billiard inside an ellipse. Thus one has Graves's theorem: {\it wrapping a closed non-stretchable string around an ellipse yields a confocal ellipse}. 

\begin{figure}[hbtp]
\centering
\includegraphics[height=1.7in]{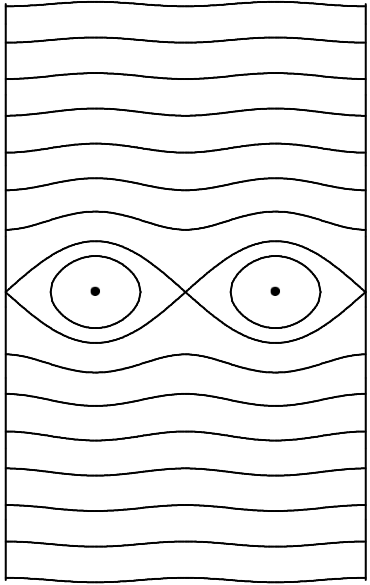}
\caption{Phase portrait of the billiard map in an ellipse}
\label{portrait}
\end{figure}

The space of rays $A$ that intersect an ellipse  is topologically a cylinder, and the billiard system inside the ellipse is an area preserving transformation $T: A\to A$. The cylinder is foliated by the invariant curves of the map $T$ consisting of the rays tangent to confocal conics, see Figure \ref{portrait}. 

The curves that go around the cylinder correspond to the rays that are tangent to confocal ellipses, and the curves that form `the eyes' to the rays that are tangent to confocal hyperbolas. A singular curve consists of the rays through the foci, and the two dots to the 2-periodic back and forth orbit along the minor axis of the ellipse.

One can choose a cyclic parameter, say, $x$ modulo 1, on each invariant circle, such that the map $T$ becomes a shift $x \mapsto x+c$, where the constant $c$ depends on the invariant curve.  Here is this construction (a particular case of the Arnold-Liouville theorem in the theory of integrable systems).

Choose a function $H$  whose level curves are the invariant curves that foliate $A$, and consider its Hamiltonian vector field sgrad $H$ with respect to the area form $\omega$. This vector field is tangent to the invariant curves, and the desired coordinate $x$ on these curves is the one in which sgrad $H$ is a constant vector field $d/dx$. Changing $H$ scales the coordinate $x$ on each invariant curve and, normalizing the `length' of the invariant curves to 1, fixes $x$ uniquely up to an additive constant. In other words, the 1-form $dx$ is well defined on each invariant curve. 

The billiard map $T$ preserves the area form and the invariant curves, therefore its restriction to each curve preserves the measure $dx$, hence, is a shift $x \mapsto x+c$.

An immediate consequence is the Poncelet Porism: if a billiard trajectory in an ellipse closes up after $n$ reflections, then 
$nc \equiv 0$ mod 1, and hence all trajectories with the same caustic close up after $n$ reflections. 

\begin{figure}[hbtp]
\centering
\includegraphics[height=1.8in]{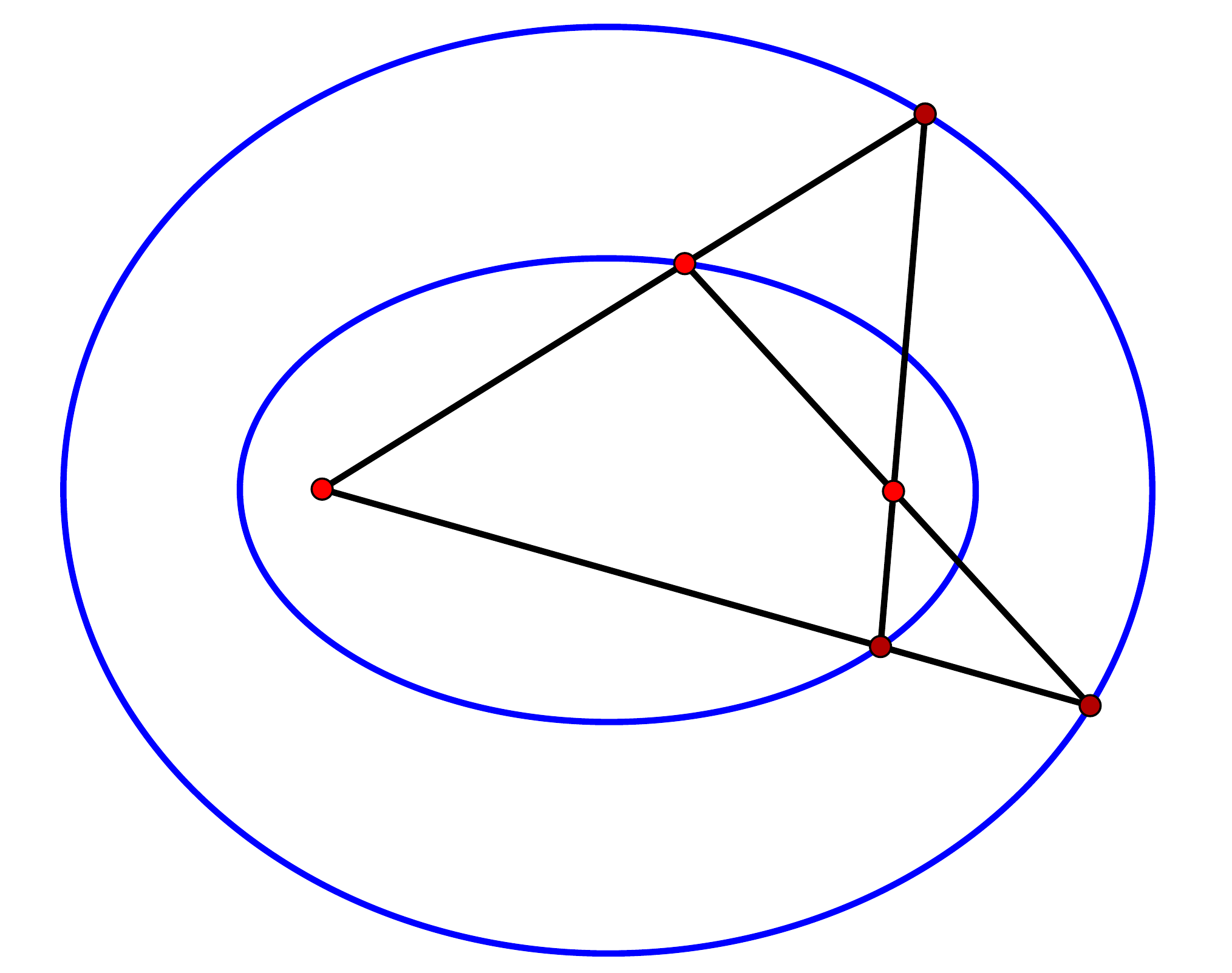}\quad 
\includegraphics[height=1.8in]{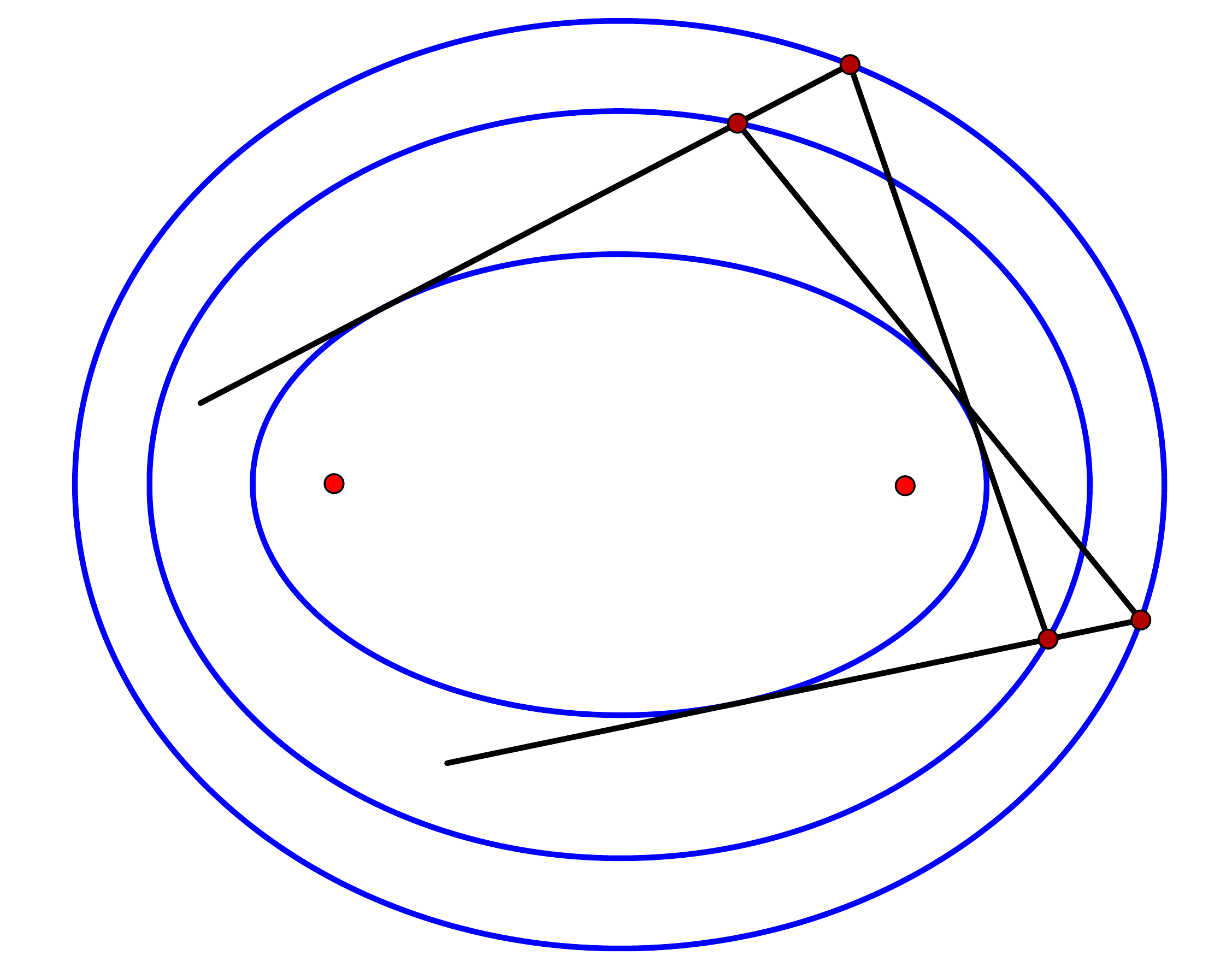}
\caption{Left: the billiard reflections of the rays from a focus 
in two confocal ellipses commute. Right: the general case.}
\label{commute}
\end{figure}

Note that the invariant measure $dx$ on the invariant curves does not depend on the choice of the billiard ellipse from 
a confocal family: the confocal ellipses share their caustics. This implies that the billiard transformations with respect to two confocal ellipses commute: restricted to a common caustic, both are shifts in the same coordinate system. This statement can be  considered as a configuration theorem; see Figure \ref{commute}.

To summarize, an ellipse is a billiard caustic for the confocal family of ellipses. It carries a coordinate $x$, defined up to an additive constant, in which the billiard reflection in confocal ellipses is given by $x \mapsto x+c$. We refer to the coordinate $x$ as the canonical coordinate.

Consider an ellipse $\gamma$, and let $x$ be the canonical  coordinate on it. Define coordinates in the exterior of the ellipse: the coordinates of a point $X$ outside of $\gamma$ are the coordinates $x_1$ and $x_2$ of the tangency points of the tangent lines from $X$ to $\gamma$. Let us call $(x_1,x_2)$ the string coordinates of point $X$. 
The confocal ellipses are given by the equations $x_2-x_1=$ const.

\begin{lemma} \label{hyp}
The confocal hyperbolas have the equations $x_2+x_1=$ const. 
\end{lemma}

\begin{figure}[hbtp]
\centering
\includegraphics[height=2in]{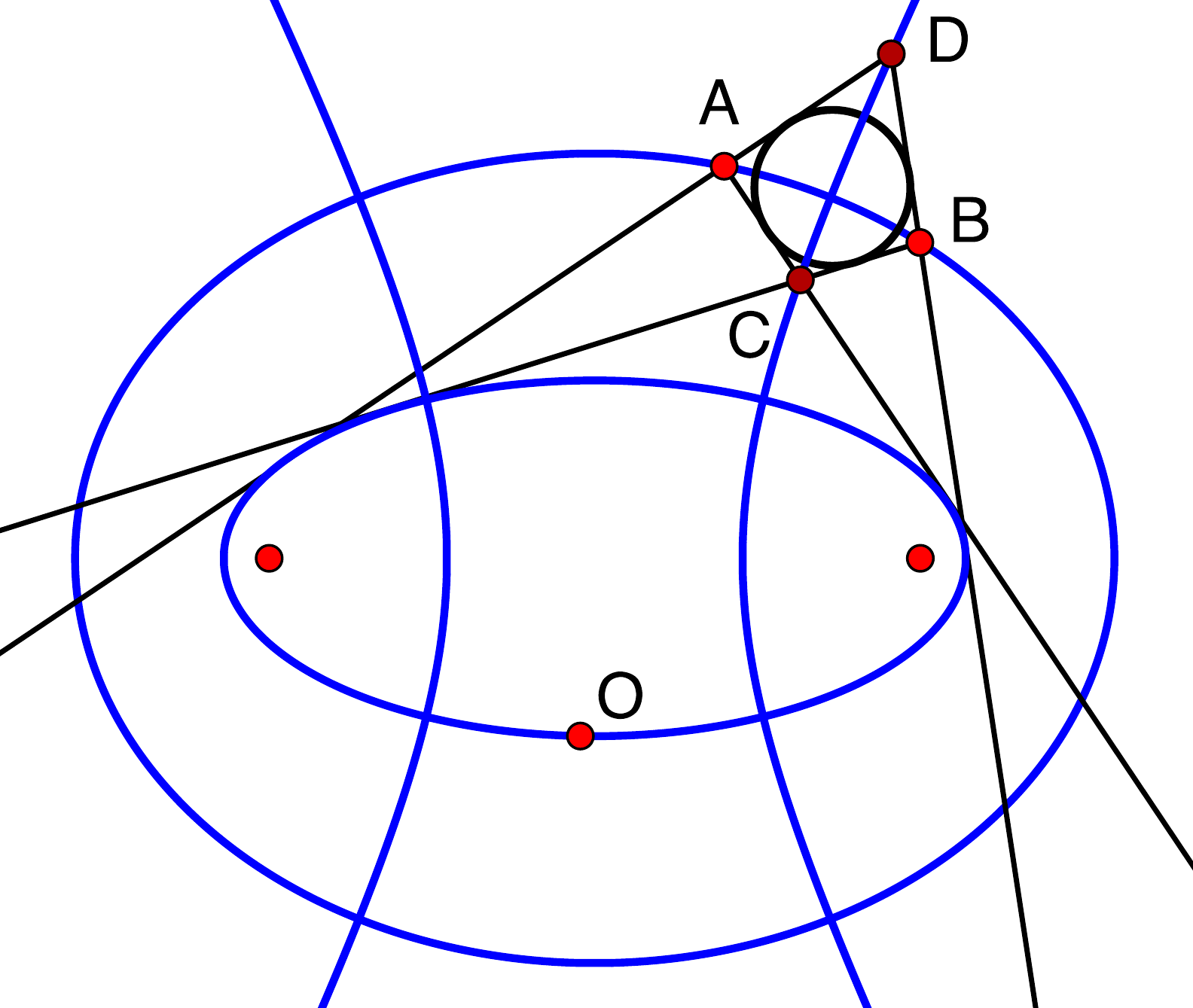}
\caption{Two pairs of tangents from an ellipse to a confocal ellipse.}
\label{circlequad}
\end{figure}

\proof Consider Figure \ref{circlequad}. Let the canonical  coordinates of the tangency points on the inner ellipse, from left to right, be $x_1,x_2,x_3,x_4$, so that the string coordinates are as follows:
$$
A(x_1,x_3),\ B(x_2,x_4),\ C(x_2,x_3),\ D(x_1,x_4).
$$
 Since $A$ and $B$ are on a confocal ellipse, $x_4-x_2=x_3-x_1$, and hence $x_2+x_3=x_1+x_4$.

By the billiard property, the arc of an ellipse $AB$ bisects the angles $CAD$ and $CBD$. Therefore, in the limit  $A\to B$, the infinitesimal quadrilateral $ABCD$ becomes a kite: the diagonal $AB$ is its axis of symmetry. Hence $AB \perp CD$, and the locus of points given by the equation $x_1+x_4=$ const and containing points $C$ and $D$  is orthogonal to the ellipse through points $A$ and $B$. Therefore this locus is a confocal hyperbola. 
\proofend

The next result is due to Reye and Chasles.

\begin{theorem} \label{inscribed}
Let $A$ and $B$ be two points on an ellipse. Consider  the quadrilateral $ABCD$, made by the pairs of tangent lines from $A$ and $B$ to a confocal ellipse. Then its other vertices, $C$ and $D$, lie on a confocal hyperbola, and the quadrilateral 
 is circumscribed about a circle,  see Figure \ref{circlequad}.
\end{theorem}

\proof In the notation of the proof of the preceding lemma, $x_2+x_3=x_1+x_4$, hence 
points $C$ and $D$ lie on a confocal hyperbola.
Furthermore, in terms of the string construction,
$$
f(A)+g(A)=f(B)+g(B),\ \ f(C)-g(C)=f(D)-g(D),
$$
hence
$$
f(D)-f(A) - g(A)+g(C)+f(B)-f(C)-g(D)+g(B)=0,
$$
or $|AD|-|AC|+|BC|-|BD|=0$.
This is necessary and sufficient for the quadrilateral $ABCD$ to be circumscribed.
\proofend

Now, consider a Poncelet $n$-gon, an $n$-periodic billiard trajectory in the ellipse $\Gamma$. Oner can choose the canonical coordinates of the tangency points of the sides of the polygon with the confocal ellipse $\gamma$   to be 
$$
 0,\ \frac{1}{n},\ \frac{2}{n},\ \dots,\ \frac{n-1}{n}.
 $$
 Then the string coordinates  of the points of the concentric set $P_k$ are
 $$
 \left(0, \frac{k}{n}\right), \left(\frac{1}{n}, \frac{k+1}{n}\right), \left(\frac{2}{n}, \frac{k+2}{n}\right), \ldots ,
 $$
 that is, their difference equals $k/n$,  a constant. It follows that $P_k$ lies on a confocal ellipse. Likewise for the radial sets $Q_k$, 
 proving the first claim of Theorem \ref{Pgrid}.
 
\begin{figure}[hbtp]
\centering
\includegraphics[height=2.5in]{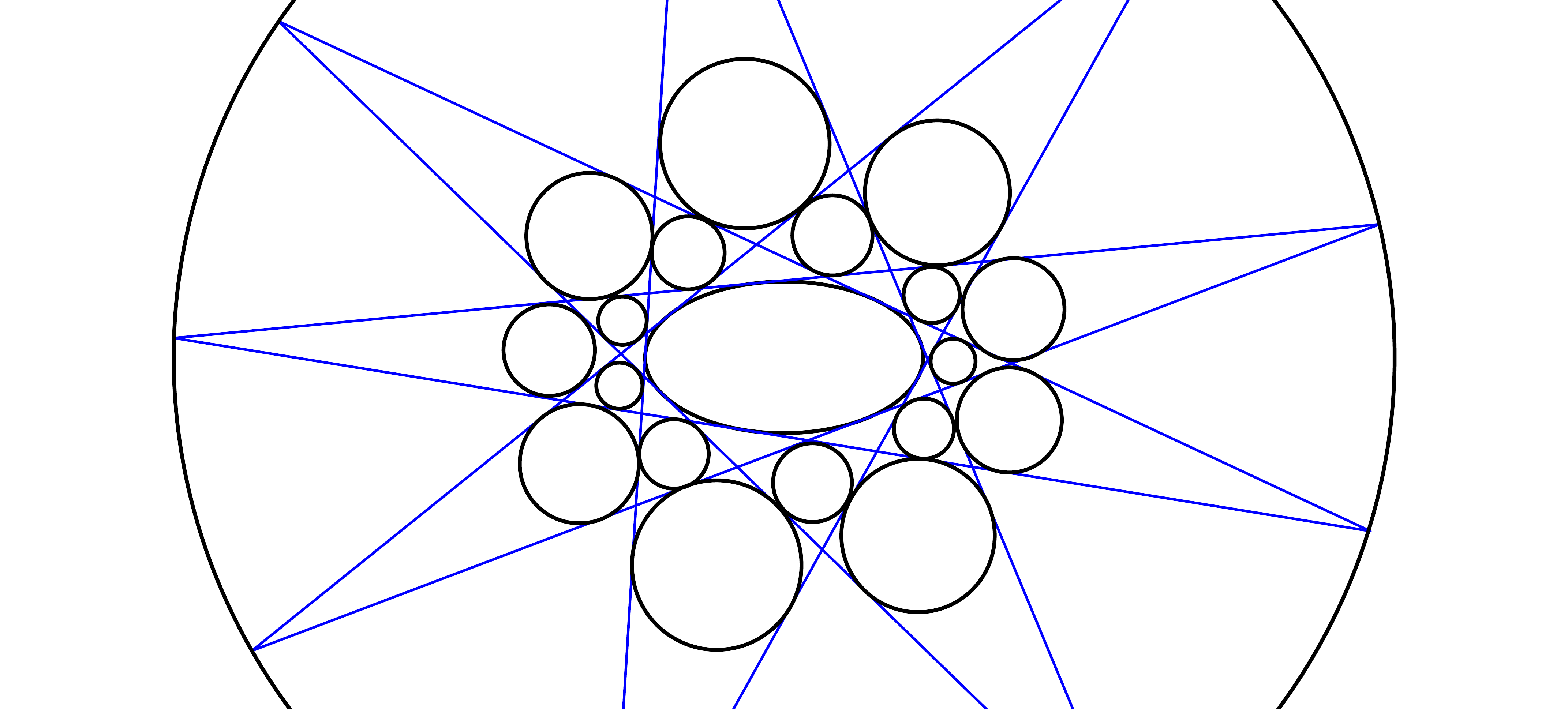}
\caption{Poncelet grid of circles.}
\label{circles}
\end{figure}

Theorem \ref{inscribed} implies that each quadrilateral of the Poncelet grid is circumscribed, see Figure \ref{circles}. We refer to \cite{AB} for circle patterns related to conics. 
 
 Next, we prove the second claim of Theorem \ref{Pgrid}.
 The confocal family of conics is given by the equation
 $$
 \frac{x_1^2}{a_1^2+\lambda}+\frac{x_2^2}{a_2^2+\lambda}=1,
 $$
 where $\lambda$ is a parameter. Its dual family is the pencil
 $$
 (a_1^2+\lambda) x_1^2+(a_2^2+\lambda) x_2^2=1
 $$
 that consists of the conics that share four points, possibly, complex. Hence the confocal family consists of the conics  that share four tangent lines, also possibly complex.

To prove the last claim of Theorem \ref{Pgrid}, we need the following classical result. Let $\gamma$ and $\Gamma$ be confocal ellipses, centered at the origin and symmetric with respect to the coordinate axes, and let $A$ be the diagonal matrix with positive entries that takes $\gamma$ to $\Gamma$.

\begin{lemma}[Ivory] \label{Ivory}
For every point $P \in \gamma$, the points $P$ and $A(P)$ lie on a confocal hyperbola.
\end{lemma}

Let us show that the linear map $A$ takes $P_k$ to $P_m$ or to its centrally symmetric set; the argument for the radial sets is similar. 

It is convenient to change the string coordinates $(x,y)$ to $u=(x+y)/2, v= (y-x)/2$. The $(u,v)$-coordinates of the points of the sets $P_k$ and $P_m$ are
$$
\left(\frac{k}{2n}+\frac{j}{n}, \frac{k}{2n}\right),\  
\left(\frac{m}{2n}+\frac{j}{n}, \frac{m}{2n}\right)\ (j=0,1,\dots,n-1).
$$
We know that $P_k$ and $P_m$ lie on confocal ellipses $\gamma$ and $\Gamma$. According to Lemma \ref{Ivory}, the map  $A$ preserves the $u$-coordinate. Therefore the coordinates of the points of the set $A (P_k)$ are
$$
\left(\frac{k}{2n}+\frac{j}{n}, \frac{m}{2n}\right)\ (j=0,1,\dots,n-1).
$$
If $m$ has the same parity as $k$, this coincides with the set $P_m$, and if the parity of $m$ is opposite to that of $k$, then this set is centrally symmetric to the set $P_m$. This completes the proof.

\section{Identities in the Lie algebras of motions} \label{Lie}

It is well known that the altitudes of a Euclidean triangle are concurrent. It is a lesser known fact that an analogous theorem holds in the spherical and hyperbolic geometries. 

In this section, we describe V. Arnold's observation \cite{Ar} that these results have interpretations as the Jacobi identity in the Lie algebras of motions of the respective geometries of constant positive or negative curvatures; see also \cite{Iv,Sk}. Following \cite{Ai,To}, we shall  discuss the relation of othr classical configuration theorems with identities in these Lie algebras.

In spherical geometry, one has the duality between points and lines that assigns the pole to an equator. There are two poles of a great circle; one can make the choice of the pole unique by considering oriented great circle, or by factorizing by the antipodal involution, that is, by replacing the sphere by the elliptic plane.

This spherical duality can be expressed in terms of the cross-product: if $A$ and $B$ are two vectors in $\R^3$ representing points in the elliptic plane, then the vector $A \times B$ represents the point dual to the line $AB$. In the following argument, we do not distinguish between points and their dual lines.

\begin{figure}[hbtp]
\centering
\includegraphics[height=2.1in]{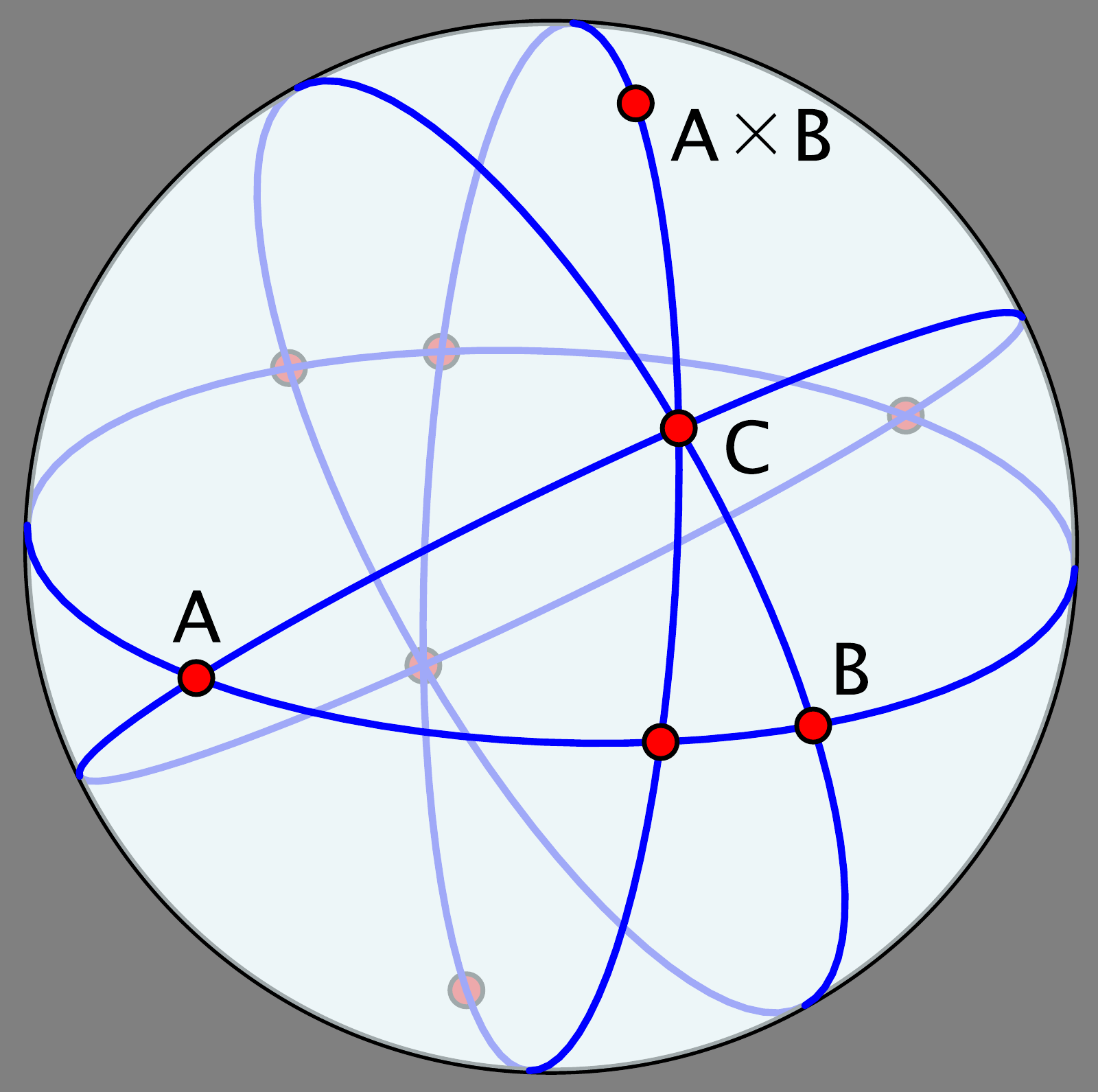}
\caption{An altitude of a spherical triangle.}
\label{sphere}
\end{figure}

Given a spherical triangle $ABC$, the altitude dropped from $C$ to $AB$ is the great circle connecting the pole of the great circle $AB$ and point $C$. Using the identification of points and lines, and cross-product, this altitude is represented by the vector $(A\times B) \times C$, see Figure \ref{sphere}. 

Two other altitudes are given by similar cross-products, and the statement that the three great circles are concurrent is equivalent to linear dependence of the these three cross-products. But
$$
(A\times B) \times C + (B\times C)\times A + (C\times A)\times B =0,
$$
the Jacobi identity for cross-product, hence the three altitudes are concurrent.

Note that the Lie algebra $(\R^3,\times)$ is $so(3)$, the algebra of motions of the unit sphere. Thus the Jacobi identity in  $so(3)$ implies the the existence of the spherical orthocenter. 

A similar, albeit somewhat more involved, argument works in the hyperbolic plane, with the Lie algebra of motions $sl(2,\R)$ replacing $so(3)$. Note that these  algebras are real forms of the  complex Lie algebra $sl(2,\C)$. 

Interestingly, the Euclidean theorem on concurrence of the three altitudes of a triangle does not seem to admit an interpretation as the Jacobi identity of the Lie algebra of motions of the plane.

Developing these ideas, T. Tomihisa \cite{To} discovered the following identity.

\begin{theorem} \label{Tid}
For every quintuple of elements of the Lie algebra $sl(2)$ (with real or complex coefficients), one has
$$
[F_1,[[F_2,F_3],[F_4,F_5]]] + [F_3,[[F_2,F_5],[F_4,F_1]]] + [F_5,[[F_2,F_1],[F_4,F_3]]] =0.
$$
\end{theorem}

Note that the indices $1,3,5$ permute cyclically, while $2$ and $4$ are frozen.

\begin{figure}[hbtp]
\centering
\includegraphics[height=1.8in]{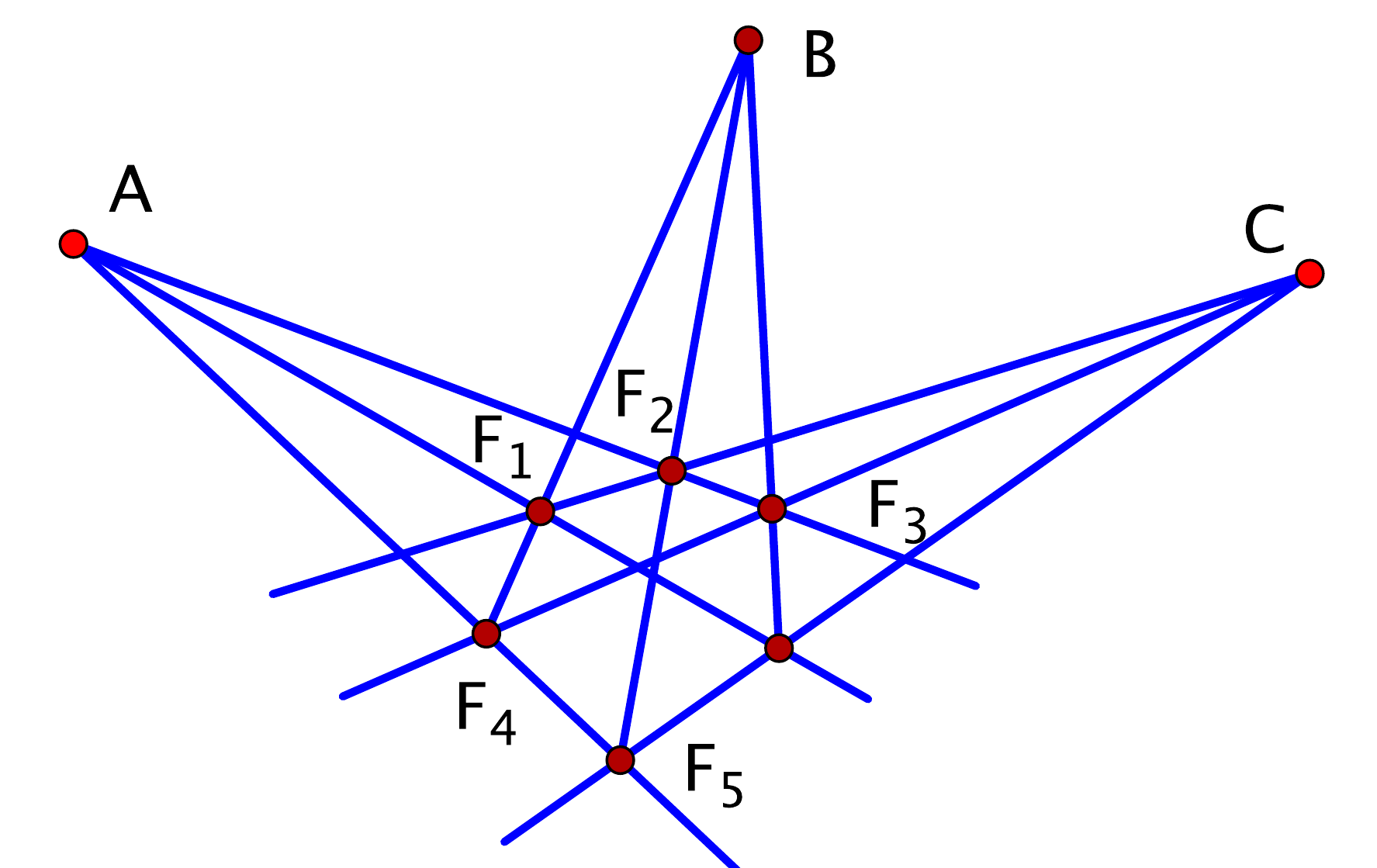}
\caption{The Tomihisa identity as the dual Pappus theorem: the lines $AF_1, BF_3$, and $CF_5$ are concurrent.}
\label{Tomihisa}
\end{figure}

As above, the Tomihisa identity can be interpreted as a configuration theorem: the Lie bracket corresponds to one of the two basic operations: connecting a pair of points by a line or intersecting a pair of lines at a point. See Figure \ref{Tomihisa} for  such an interpretation.

\section{Skewers} \label{skewers}

This section is based upon the recent paper \cite{Tab}. The main idea is that planar projective configuration theorems have space analogs where points and lines in the projective plane are replaced by lines in space, and the two operations, connecting two points by a line and intersecting two lines at a point, are replaced by taking the common perpendicular of two lines.

The {\it skewer} of two lines in 3-dimensional space is their common perpendicular.
We denote the skewer of lines $a$ and $b$ by $S(a,b)$. In Euclidean and hyperbolic spaces, a generic pair of lines has a unique skewer; in the spherical geometry, a generic pair of lines (great circles) has two skewers, similarly to a great circle on $S^2$ having two poles. We always assume that the lines involved in the formulations of the theorems are in general position.

Here is the `skewer translation' of the Pappus theorem, as depicted in Figure \ref{Pappusfig}:

\begin{theorem} \label{skPappus}
Let $a_1,a_2,a_3$ be a triple of lines with a common skewer, and let $b_1,b_2,b_3$ be another  triple of lines with a common skewer. Then the  lines 
$$
S(S(a_1,b_2),S(a_2,b_1)),\ S(S(a_1,b_3),S(a_3,b_1)),\  {\rm and}\ \ S(S(a_2,b_3),S(a_3,b_2))
$$
 share a skewer.
\end{theorem} 

This theorem, as well as in the following ones, holds in $\R^3, S^3$ and $H^3$. 

And here is the skewer version of the Desargues theorem, as depicted in Figure \ref{Desarguesfig}:

\begin{theorem}\label{skDesargues}
Let $a_1,a_2,a_3$ and $b_1,b_2,b_3$ be two triples of lines such that the lines $S(a_1,b_1), S(a_2,b_2)$ and $S(a_3,b_3)$ share a skewer. Then the lines 
$$
S(S(a_1,a_2),S(b_1,b_2)),\ S(S(a_1,a_3),S(b_1,b_3)),\  {\rm and}\ \ S(S(a_2,a_3),S(b_2,b_3))
$$
also share a skewer.
\end{theorem}

The `rules of translation' should be clear from these examples.

As a third example, consider a configuration theorem that involves polarity, namely, the theorem that the three altitudes of a triangle are concurrent that was discussed in Section \ref{Lie}. In its skewer version, one does not distinguish between polar dual objects, such as a great circle and its pole. This yields

\begin{theorem} \label{skAlt}
Given three lines $a,b,c$, the  lines
$$
S(S(a,b),c),\ S(S(b,c),a),\ \ {\rm and}\ \ S(S(c,a),b)
$$
share a skewer.
\end{theorem}

\begin{figure}[hbtp]
\centering
\includegraphics[height=1.8in]{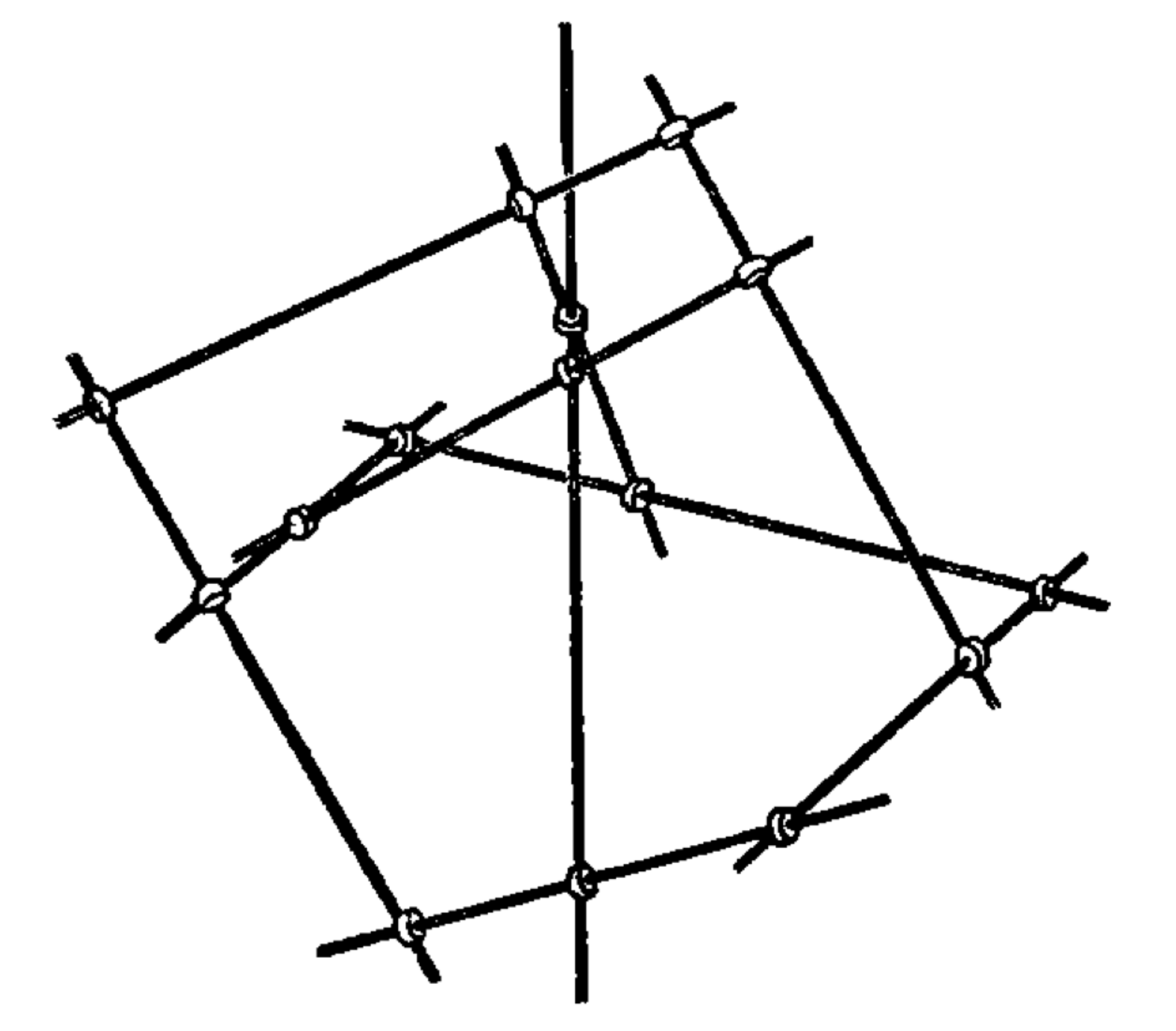}
\caption{Petersen-Morley configuration of ten lines.}
\label{ten}
\end{figure}

This is the Petersen-Morley, also known as Hjelmslev-Morley, theorem \cite{Mo}. An equivalent formulation:
 {\it the common normals of the opposite sides of a rectangular hexagon have a common normal}. See Figure \ref{ten}, borrowed from \cite{Mo2}.

Denote the 2-parameter family of lines that meet a given line $\ell$ at right angle by ${\cal N}_{\ell}$. The sets ${\cal N}_{\ell}$ plays the role of lines in the skewer versions of configuration theorems. Two-parameter families of lines in 3-space are called congruences.

Next we describe line analogs of circles. Let $\ell$  be an oriented  line in 3-space (elliptic, Euclidean, or hyperbolic). Let $G_\ell$ be the 2-dimensional subgroup of the group of orientation preserving isometries that preserve $\ell$. The orbit $G_\ell(m)$ of an oriented line $m$ is called the {\it axial congruence} with $\ell$ as axis (an analog of the center of a circle). 

In $\R^3$, the lines of an axial congruence with axis $\ell$ are at equal distances from $\ell$ and make equal angles with it. In the hyperbolic space, one can define the so-called complex distance between oriented lines, see \cite{Ma}. The complex distance between the lines of an axial congruence and its axis is constant.

Axial congruences share the basic properties of circles: if two generic axial congruences  share a line, then they share a unique other line; and three generic oriented lines belong to a unique axial congruence.

The next result is a skewer analog of the Pascal theorem, see Figure \ref{Pascalfig}, in the particular case when the conic is a circle.

\begin{theorem} \label{skPascal}
Let $A_1,\ldots,A_6$ be lines from an axial congruence. Then
$$
S(S(A_1,A_2),S(A_4,A_5)),\ S(S(A_2,A_3),S(A_5,A_6)),    \  {\rm and}\  S(S(A_3,A_4),S(A_6,A_1))
$$  
share a skewer.
\end{theorem}

As another, lesser known, example, consider the  Clifford's Chain of Circles. This chain of theorems starts with a collection of concurrent circles labelled $1,2,3,\ldots, n$.  The intersection point of the circles $i$ and $j$ is labelled $ij$. The circle through points $ij, jk$ and $ki$ is labelled $ijk$. 

\begin{figure}[hbtp]
\centering
\includegraphics[height=3.2in]{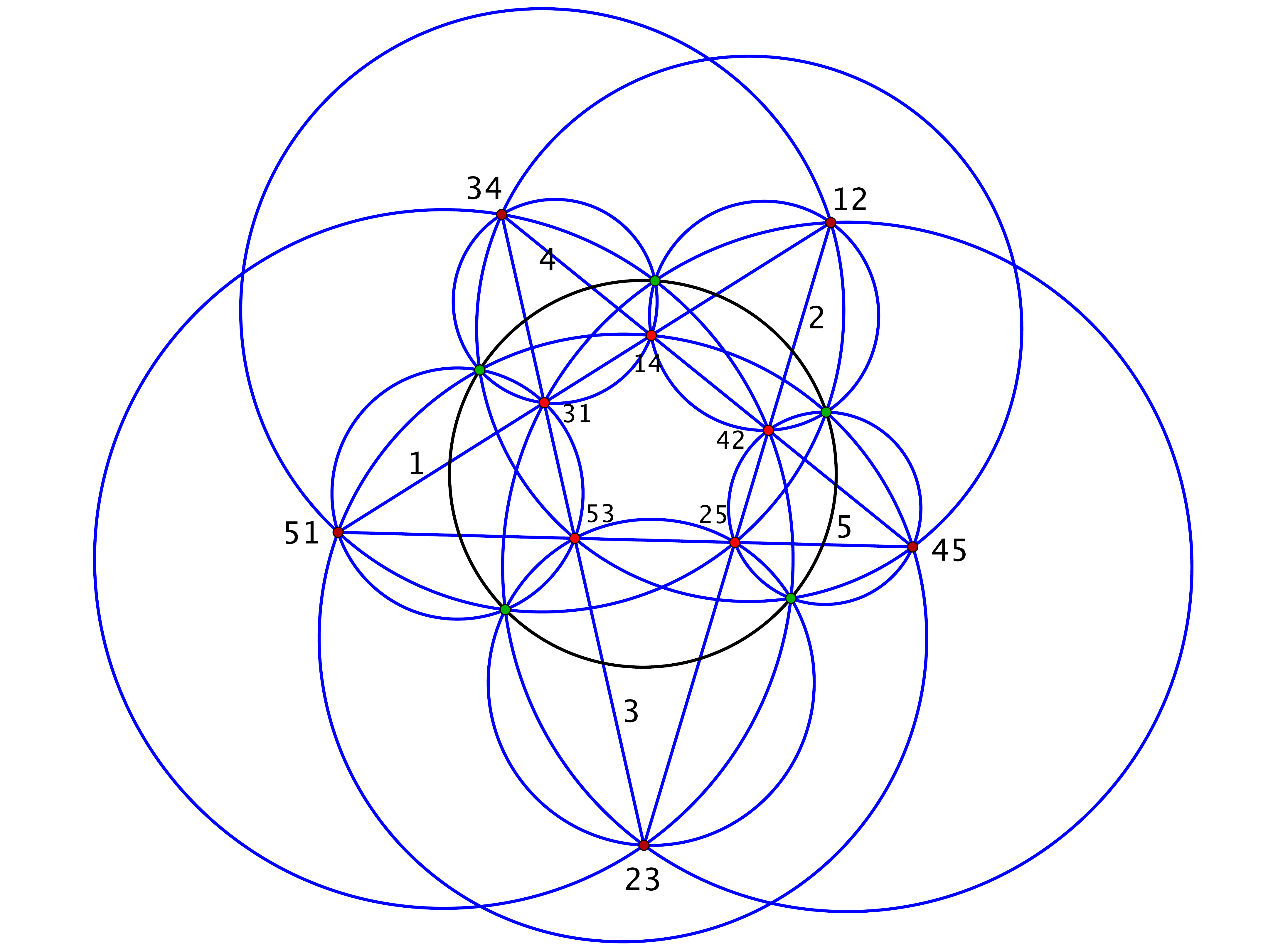}
\caption{Clifford's Chain of Circles ($n=5$).}
\label{Clifford}
\end{figure}

The first statement of the theorem is that the circles $ijk, jkl, kli$ and $lij$ share a point;  this point is labelled $ijkl$. The next statement is that the points $ijkl, jklm, klmi, lmij$ and $mijk$ are concyclic; this circle is labelled $ijklm$. And so on, with the alternating claims of being concurrent and concyclic; see \cite{Coo,Mo2}, and Figure \ref{Clifford} where the initial circles are represented by lines (circles of infinite radius sharing a point at infinity).

The next theorem, in the case of  $\R^3$, is due to Richmond \cite{Ri}. 

\begin{theorem} \label{skClifford}
1) Consider axial congruences ${\cal C}_i,\ i=1,2,3,4$,  sharing a line.  For each pair of indices $i,j \in \{1,2,3,4\}$, denote by $\ell_{ij}$ the  line shared by ${\cal C}_i$ and ${\cal C}_j$. For each triple of indices $i,j,k \in \{1,2,3,4\}$, denote by ${\cal C}_{ijk}$  the axial congruence containing the lines $\ell_{ij},\ell_{jk},\ell_{ki}$. Then the congruences ${\cal C}_{123}, {\cal C}_{234}, {\cal C}_{341}$ and ${\cal C}_{412}$ share a line. \\
2) Consider axial congruences ${\cal C}_i,\ i=1,2,3,4,5$,  sharing a line. Each four of the indices determine a line, as described in the previous statement of the theorem. One obtains five lines, and they all belong to an axial congruence.\\
3) Consider axial congruences ${\cal C}_i,\ i=1,2,3,4,5,6$,  sharing a line. Each five of them determine an axial congruence, as described in the previous statement of the theorem. One obtains six axial congruences, and they all share a line.
And so on...
\end{theorem} 

Next one would like to define  line analogs of conics. A first step in this direction is made in \cite{Tab}, but much more work is needed. In particular, one would like to have skewer analogs of various configuration theorems involving conics, including the Pascal theorem and the whole hexagrammum mysticum, the 
Poncelet Porism, and the theorems described in Section \ref{pentalike}. As of now, this is an open problem.

Now we outline two approaches to proofs of the above theorems and the skewer versions of other planar configuration theorems. The first approach is by way of the  spherical geometry, and the second via the hyperbolic geometry. Either approach implies the results in all three classical geometries by `analytic continuation'.
This analytic continuation principle is well known in geometry; see, e.g.,  \cite{AP1,Pa} where it is discussed in detail.

\paragraph{Elliptic approach.} The space of oriented great circles in $S^3$, or lines in the elliptic space $\RP^3$, is he  Grassmannian $G(2,4)$ of oriented 2-dimensional subspaces in $\R^4$. Below we collect pertinent facts concerning this Grassmannian.

To every oriented line $\ell$ in $\RP^3$ there corresponds its dual oriented line $\ell^*$: the respective oriented planes in $\R^4$ are the orthogonal complements of each other. The dual lines are equidistant and they have infinitely many skewers.

The Grassmannian is a product of two spheres:  $G(2,4)=S^2_-\times S^2_+$. This provides an identification of an oriented line in $\RP^3$ with a pair of points of the unit sphere $S^2$: $\ell\leftrightarrow (\ell_-,\ell_+)$. 
The antipodal involutions of the spheres $S^2_-$ and $S^2_+$ generate the action of the Klein group $\Z_2\times\Z_2$ on the space of oriented lines generated by reversing the orientation of a line and by taking the dual line.

Two lines  $\ell$ and $m$ intersect at  right angle if and only if $d(\ell_-,m_-)=d(\ell_+,m_+)=\pi/2$, where $d$ denotes the spherical distance in $S^2$. A line $n$ is a skewer of lines $\ell$ and $m$ if and only if $n_-$ is a pole of the great circle $\ell_- m_-$, and $n_+$ is a pole of the great circle $\ell_+ m_+$.

The set of lines that intersect $\ell$ at  right angle coincides with the set of lines that intersect $\ell$ and $\ell^*$.
A generic pair of  lines has exactly two skewers (four, if orientation is taken into account), and they are dual to each other.

It follows that a configuration involving lines in elliptic space and their skewers can be identified with a pair of configurations on the spheres $S^2_-$ and $S^2_+$. Under this identification, the great circles of these spheres are not distinguished from their poles, just like in the proof described in Section \ref{Lie}. That is, the operation of taking the skewer of two lines is represented, on both spheres, by the cross-product. 

In this way, a configuration of lines in space becomes the direct product of the corresponding planar configurations. For example, the Petersen-Morley Theorem \ref{skAlt} splits into two statements that the altitudes of  triangles, on the spheres $S^2_-$ and $S^2_+$, are concurrent. 

\paragraph{Hyperbolic approach.} In a nutshell, a skewer configuration theorem in 3-dimensional hyperbolic space  is a complexification of a configuration theorem in the hyperbolic plane. We follow the ideas of F. Morley \cite{Mo1,Mo2}, Coxeter \cite{Co}, and V. Arnold \cite{Ar}.

Consider the hyperbolic space  in the upper halfspace model. The isometry group is $SL(2,\C)$, and the sphere at infinity (the celestial sphere of \cite{Mo1}) is the Riemann sphere $\CP^1$. 

A line in $H^3$ intersects the sphere at infinity at two points, hence the space of (non-oriented) lines is the configuration space of unordered pairs of points. As we mentioned in Section \ref{itSteiner}, $S^2(\CP^1)=\CP^2$, namely, 
to a pair of points in the projective line one assigns the binary quadratic form having zeros at these points:
$$
(a_1:b_1,a_2:b_2) \longmapsto (a_1y-b_1x)(a_2y-b_2x).
$$
Thus a line in $H^3$ can be though of as a complex binary quadratic form, up to a  factor. 

The space of binary quadratic forms $ax^2+2bxy+cy^2$ has the discriminant quadratic form $\Delta=ac-b^2$ and the respective bilinear form. The equation $\Delta=0$ defines the diagonal of $S^2(\CP^1)$; this is a conic in $\CP^2$ that does not correspond to lines in $H^3$.

 The next result is contained in \S 52 of \cite{Mo2}.

\begin{lemma} \label{Jacobian}
Two lines in $H^3$ intersect at  right angle if and only if the respective binary quadratic forms 
$f_i=a_i x^2 + 2 b_i xy + c_i y^2,\ i=1,2$, are orthogonal with respect to $\Delta$:
\begin{equation} \label{ort}
a_1c_2-2b_1b_2+a_2c_1=0.
\end{equation}
If two lines correspond to binary quadratic forms $f_i=a_i x^2 + 2 b_i xy + c_i y^2,\ i=1,2$, then their
skewer  corresponds to the Poisson bracket (the Jacobian)
$$
\{f_1,f_2\} = (a_1b_2-a_2b_1)x^2 + (a_1c_2-a_2c_1) xy + (b_1c_2-b_2c_1) y^2.
$$
\end{lemma} 

If $(a_1:b_1:c_1)$ and $(a_2:b_2:c_2)$ are homogeneous coordinates in the projective plane and the dual projective plane, then (\ref{ort}) describes the incidence relation between points and lines. In particular, the set of lines in $H^3$ that meet a fixed line at  right angle corresponds to a line in $\CP^2$.

Suppose a configuration theorem involving polarity is given in $\RP^2$. The projective plane with a conic provides the projective model of the hyperbolic plane, so the configuration is realized in $H^2$. Consider the complexification, the respective configuration theorem in $\CP^2$ with the polarity induced by $\Delta$. According to Lemma \ref{Jacobian}, this yields  a configuration of lines in $H^3$ such that the pairs of incident points and lines correspond to pairs of lines intersecting at right angle. 

\begin{remark}[On Lie algebras] \label{Liealgebraic}
{\rm From the point of view of the identities in Lie algebras, discussed in Section \ref{Lie}, the relation between configuration theorems in the hyperbolic plane and the hyperbolic space is the relation between $sl(2,\R)$ and $sl(2,\C)$: an identity in the former implies the same identity in the latter.

As to the Lie algebras in space, in the elliptic case, the Lie algebra of motions is $so(4)=so(3)\oplus so(3)$, and in the hyperbolic case, it is $sl(2,\C)$. Accordingly, an elliptic skewer configuration splits into two configurations in $S^2$, and a hyperbolic skewer configuration is obtained from a configuration in $H^2$ by complexification.
}
\end{remark}

We finish the section by discussing two results concerning lines in 3-space that do not follow the above described general pattern. The first of them is the skewer version of the Sylvester Problem.  

Given a finite set $S$ of points in the plane, assume that the line through every pair of points in $S$ contains at least one other point of $S$. J.J. Sylvester asked in 1893 whether $S$ necessarily consists of collinear points. See \cite{BM} for the history of this problem and its generalizations.

In $\RP^2$, the Sylvester Problem, along with its dual, has an affirmative answer (the Sylvester-Galai theorem), but in $\CP^2$ one has a counter-example: the 9 inflection points of a cubic curve (of which at most three can be real, according to a theorem of Klein), connected by 12 lines.

The skewer version of the Sylvester Problem concerns a finite collection of pairwise skew lines in space such that the skewer of any pair intersects at least one other line at right angle. The question is whether a collection of lines with this skewer Sylvester property necessarily consists of the lines that intersect some line at right angle. 

\begin{theorem} \label{Sylvth}
The skewer version of the Sylvester-Galai theorem holds in the elliptic and Euclidean geometries, but fails in the hyperbolic geometry.
\end{theorem}

\proof In the elliptic case, we argue as in the above described elliptic proof. A collection of lines becomes two collections of points,  in $\RP^2_-$ and in $\RP^2_+$, and the skewer Sylvester property implies that each of these sets has the property that the line through every pair of points contains another point, so one applies the Sylvester-Galai theorem on each sphere.

In the hyperbolic case, we argue as in the hyperbolic proof. Let $a_1,\ldots,a_9$ be the nine inflection points of a cubic curve in $\CP^2$, and let $b_1,\ldots,b_{12}$ be the respective lines. Let $b_1^*,\ldots,b_{12}^*$ be the polar dual points. Then the points $a_i$ correspond to nine lines in $H^3$, and the points $b_j^*$ to their skewers. We obtain a collection of nine lines that has the skewer Sylvester property but does not possess a common skewer.

In the intermediate case of $\R^3$, the argument is due to V. Timorin (private communication).

 Let us add to $\R^3$ the plane at infinity $H$; the points of $H$ are the directions of lines in space. One has a polarity in $H$ that assigns to a direction the set of the orthogonal directions, a line in $H$.

Therefore, if three lines in $\R^3$ share a skewer, then their intersections with the plane $H$ are collinear. 
Let $L_1,\ldots, L_n$ be a collection of  lines with the skewer Sylvester property. Then, by the Sylvester-Galai theorem in $H$, the points $L_1 \cap H,\ldots, L_n \cap H$  are collinear. This means that the lines $L_1,\ldots, L_n$ lie in parallel planes, say, the horizontal ones. 

Consider the vertical projection of these lines. We obtain a finite collection of non-parallel lines in the plane such that through the intersection point of any two there passes at least one other line. By the dual Sylvester-Galai theorem, all these lines are concurrent. Therefore the respective horizontal lines in $\R^3$ share a vertical skewer. 
\proofend

The second result is a different skewer version of the Pappus theorem.

\begin{theorem} \label{othPapp}
Let $\ell$ and $m$ be a pair of skew lines. Choose a triple of points $A_1,A_2,A_3$ on $\ell$ and a triple of points $B_1,B_2,B_3$ on $m$. Then the lines 
$$
S((A_1 B_2), (A_2 B_1)), \  S((A_2 B_3), (A_3 B_2)),\  {\rm and}\ \ S((A_3 B_1), (A_1 B_3))
$$
share a skewer.
\end{theorem}

We proved this result, in the hyperbolic case, by a brute force calculation using the approach to hyperbolic geometry, developed in \cite{Fe}; see \cite{Tab} for details. It is not clear whether this theorem is a part of a general pattern.

Let us close with  inviting the reader to mull over the skewer versions of other constructions of planar projective geometry. For example, one can define the skewer pentagram map that acts on cyclically ordered tuples of lines in space:
$$
\{L_1, L_2,\ldots \} \mapsto \{S(S(L_1,L_3),S(L_2,L_4)), S(S(L_2,L_4),S(L_3,L_5)), \ldots \}
$$
Is this map completely integrable?


\begin{thebibliography}{99}

\bibitem{AP1} N. A'Campo, A. Papadopoulos. {\it Transitional geometry.} in {\it Sophus Lie and Felix Klein: The Erlangen Program and its Impact in Mathematics and Physics}, 217--235. European Math. Soc., Z\"urich, 2015.

\bibitem{Ai} F. Aicardi. {\it Projective geometry from Poisson algebras.} J. Geom. Phys. {\bf 61} (2011), 1574--1586.

\bibitem{AB} A. Akopyan, A. Bobenko. {\it Incircular nets and confocal conics}. arXiv arXiv:1602.04637.

\bibitem{Ar} V. Arnold. {\it Lobachevsky triangle altitude theorem as the Jacobi identity in the Lie algebra of quadratic forms on symplectic plane}. J. Geom. Phys. {\bf 53} (2005), 421--427.

\bibitem{Be1} M. Berger. {\it Geometry. I. II.}  Springer-Verlag, Berlin, 1987. 

\bibitem{Be} M. Berger. {\it Geometry revealed.  A Jacob's ladder to modern higher geometry.} Springer, Heidelberg, 2010.

\bibitem{BS} A. Bobenko, Yu. Suris. {\it Discrete differential geometry. 
Integrable structure.} Amer. Math. Soc., Providence, RI, 2008. 

\bibitem{BM} P. Borwein, W. O. J. Moser. {\it A survey of Sylvester's problem and its generalizations.} Aequationes Math. {\bf 40} (1990), 111--135. 

\bibitem{CR1} J. Conway, A. Ryba. {\it The Pascal mysticum demystified.} Math. Intelligencer {\bf 34} (2012), no. 3, 4--8. 

\bibitem{CR2} J. Conway, A. Ryba. {\it Extending the Pascal mysticum.} Math. Intelligencer {\bf 35} (2013), no. 2, 44--51.

\bibitem{Coo} J. L.  Coolidge. {\it A treatise on the circle and the sphere.} Chelsea Publ. Co., Bronx, N.Y., 1971.

\bibitem{Co} H. S. M. Coxeter. {\it The inversive plane and hyperbolic space.}
 Abh. Math. Sem. Univ. Hamburg {\bf 29} (1966), 217--242.
 
 \bibitem{DR} V.  Dragovi\'c, M.  Radnovi\'c. {\it Poncelet porisms and beyond. 
Integrable billiards, hyperelliptic Jacobians and pencils of quadrics.}  Birkh\"auser/Springer, Basel, 2011.

\bibitem {Fe} W. Fenchel. {\it Elementary geometry in hyperbolic space.} Walter de Gruyter, Berlin, 1989.

\bibitem {Fl} L. Flatto. {\it Poncelet's theorem.} Amer. Math. Soc., Providence, RI, 2009.

\bibitem{FT} D. Fuchs, S. Tabachnikov. {\it Mathematical omnibus. Thirty lectures on classic mathematics.} Amer. Math. Soc., Providence, RI, 2007. 

\bibitem{GSTV} M. Gekhtman, M. Shapiro, S. Tabachnikov, and A. Vainshtein, {\it Integrable cluster dynamics of directed networks and pentagram maps}, with appendix by A. Izosimov, Adv. Math., in print.

\bibitem{Gl} M. Glick. {\it The pentagram map and $Y$-patterns.} Adv. Math. {\bf 227} (2011), 1019--1045.

\bibitem{GP} M. Glick, P. Pylyavskyy. {\it $Y$-meshes and generalized pentagram maps.} Proc. Lond. Math. Soc.  {\bf 112} (2016),  753--797.

\bibitem {Gr} B. Gr\"unbaum. {\it Configurations of points and lines}.  Amer. Math. Soc., Providence, RI, 2009.

\bibitem{Ha} T. Hatase. {\it Algebraic Pappus Curves.}  
Ph.D. Thesis, Oregon State University, 2011.

\bibitem{HC} D. Hilbert, S. Cohn-Vossen. {\it Geometry and the imagination.} Chelsea Publishing Co, New York, N. Y., 1952.

\bibitem {Ho} P. Hooper. {\it From Pappus' theorem to the twisted cubic.} Geom. Dedicata {\bf 110} (2005), 103--134. 

\bibitem{Iv} N. Ivanov. {\it Arnol'd, the Jacobi identity, and orthocenters.}
Amer. Math. Monthly {\bf 118} (2011),  41--65.


\bibitem{Ki} F. Kissler. {\it A family of representations for the modular group}. Master Thesis, Heidelberg, 2016.

\bibitem{KT} V. Kozlov and D. Treshchev. {\it Billiards. A Genetic Introduction to the Dynamics of Systems with Impacts.} Amer. Math. Soc., Providence, RI, 1991.

\bibitem{LT} M. Levi, S.  Tabachnikov. {\it The Poncelet grid and billiards in ellipses}. Amer. Math. Monthly {\bf 114} (2007), 895--908.

\bibitem {Ma} A. Marden. {\it Outer circles. An introduction to hyperbolic 3-manifolds.} Cambridge Univ. Press, Cambridge, 2007. 

\bibitem{Mo} F. Morley. {\it On a regular rectangular configuration of ten lines}. Proc. London Math. Soc. s1-29 (1897), 670--673.

\bibitem{Mo1} F. Morley. {\it The Celestial Sphere.} Amer. J. Math. {\bf 54} (1932),  276--278.

\bibitem{Mo2} F. Morley, F. V. Morley. {\it Inversive geometry}. G. Bell \& Sons, London, 1933.

\bibitem{OST1} V. Ovsienko, R. Schwartz,  S. Tabachnikov, {\it
    The pentagram map: a discrete integrable system,} Comm. Math.
  Phys. {\bf 299} (2010), 409--446.
  
\bibitem{OST2} V. Ovsienko, R. Schwartz, and S. Tabachnikov, {\it Liouville-Arnold integrability of the pentagram map on closed polygons}, Duke Math. J. {\bf 162} (2013), 2149--2196.

\bibitem{Pa} {\it Strasbourg master class on geometry.} 
A. Papadopoulos ed.  European Math. Soc., Z\"urich, 2012.

\bibitem{RG} J. Richter-Gebert. {\it Perspectives on projective geometry. A guided tour through real and complex geometry.} Springer, Heidelberg, 2011.

\bibitem{Ri} J. F. Rigby. {\it Pappus lines and Leisenring lines.} J. Geom. {\bf 21} (1983), 108--117. 

\bibitem {Sch92} R. Schwartz. {\it The pentagram map}. Experiment. Math. {\bf 1} (1992), 71--81.

\bibitem {Sch93} R. Schwartz. {\it Pappus' theorem and the modular group}. Inst. Hautes \'Etudes Sci. Publ. Math. {\bf 78} (1993), 187--206 (1994).

\bibitem {Sch01} R. Schwartz. {\it Desargues theorem, dynamics, and hyperplane arrangements}. Geom. Dedicata {\bf 87} (2001), 261--283. 

\bibitem{Sch07} R. Schwartz. {\it The Poncelet grid.} Adv. Geom. {\bf 7} (2007), 157--175. 

\bibitem{Sch08} R. Schwartz. {\it Discrete monodromy, pentagrams, and the method of condensation.} J. Fixed Point Theory Appl. {\bf 3} (2008),   379--409.

\bibitem {ST1} R. Schwartz, S. Tabachnikov. {\it Elementary surprises in projective geometry}. Math. Intelligencer {\bf 32} (2010), no. 3, 31--34. 


\bibitem {Sk} M. Skopenkov. {\it Theorem about the altitudes of a triangle and the Jacobi identity} (in Russian). Matem. Prosv., Ser. 3 {\bf 11} (2007), 79--89.

\bibitem{So} F. Soloviev, {\it Integrability of the pentagram map},
  Duke Math.J. {\bf 162} (2013), 2815--2853.
  
\bibitem{Tab95} S. Tabachnikov. {\it Billiards.} Panor. Synth. No. 1, SMF, 1995. 
 
\bibitem{Tab05} S. Tabachnikov. {\it Geometry and billiards.}  Amer. Math. Soc., Providence, RI, 2005.

\bibitem{Tab} S. Tabachnikov. {\it Skewers.} Arnold Math. J. {\bf 2} (2016), 171--193.

\bibitem{To} T. Tomihisa. {\it Geometry of projective plane and Poisson structure.} J. Geom. Phys. {\bf 59} (2009),   673--684.

\bibitem{Cin} \url{http://www.cinderella.de/tiki-index.php}.

\bibitem{App} \url{https://www.math.brown.edu/~res/Java/App33/test1.html}.

\bibitem{App1} \url{https://www.math.brown.edu/~res/Java/Special/Main.html}.

\end{thebibliography}
\end{document}